\documentclass[11pt,reqno]{amsproc}
\allowdisplaybreaks
\usepackage{cite} 
\usepackage[numbers,sort&compress]{natbib}

\usepackage{makecell}
\usepackage{multirow}
\usepackage{amssymb}
\usepackage{stmaryrd}
\usepackage{amsmath}
\usepackage{color}
\usepackage{colortbl}
\usepackage{fullpage}
\usepackage{comment}
\usepackage{graphicx}
\usepackage{textcomp}
\usepackage{lineno}
\usepackage{subcaption}
\usepackage{enumerate}
\usepackage{longtable} 
\usepackage[debug=false, colorlinks=true, pdfstartview=FitV,
linkcolor=blue, citecolor=blue, urlcolor=blue]{hyperref}

\usepackage[usenames,dvipsnames,table]{xcolor}
\usepackage[most]{tcolorbox}

\usepackage[doipre={DOI:~}]{uri}
\usepackage{mathtools}
\usepackage{dirtytalk}

\usepackage{morefloats}
\usepackage{soul}
\usepackage[normalem]{ulem}
\usepackage{cancel}
\newlength{\drop}
\definecolor{amethyst}{rgb}{0.6, 0.4, 0.8}
\definecolor{burgundy}{rgb}{0.5, 0.0, 0.13}

\title{\textbf{CoolPINNs:~A Physics-informed Neural 
Network Modeling of Active Cooling in Vascular Systems}}

\author{$\textbf{Nimish~V.~Jagtap}^{1}$, 
$\textbf{M.~K.~Mudunuru}^{2}$, and $\textbf{K.~B.~Nakshatrala}^{1}$ \\
  {\small $^{1}$Department of Civil and Environmental Engineering, 
  University of Houston, Houston, Texas 77204, USA. \\
  $^{2}$Pacific Northwest National Laboratory, 
  Richland, Washington, 99352. 
  } \\
  Corresponding author: knakshatrala@uh.edu, 
  Co-corresponding author: maruti@pnnl.gov
  }


\begin{document}

\begin{titlepage}
  \drop=0.1\textheight
  \centering
  \vspace*{\baselineskip}
  \rule{\textwidth}{1.6pt}\vspace*{-\baselineskip}\vspace*{2pt}
  \rule{\textwidth}{0.4pt}\\[\baselineskip]
       {\Large \textbf{\color{burgundy}
       CoolPINNs:~A Physics-informed Neural Network Modeling of Active Cooling in Vascular Systems}}\\[0.3\baselineskip]
       \rule{\textwidth}{0.4pt}\vspace*{-\baselineskip}\vspace{3.2pt}
       \rule{\textwidth}{1.6pt}\\[0.5\baselineskip]
       \scshape
       An e-print of the paper is available on arXiv. \par
       \vspace*{0.1\baselineskip}
       Authored by \\[0.1\baselineskip]
       
        {\Large Nimish~V.~Jagtap\par}
  {\itshape Doctoral Student, Dept. of Mechanical 
  Engineering, University of Houston, Texas 77204.}\\[0.1\baselineskip]

  {\Large M.~K.~Mudunuru\par}
  {\itshape Earth Scientist, 
  Earth System Science Division, \\
  Pacific Northwest National Laboratory, 
  Richland, Washington, 99352. \\
  \textbf{phone:} +1-509-375-6645, 
  \textbf{e-mail:} maruti@pnnl.gov}\\[0.1\baselineskip]

  {\Large K.~B.~Nakshatrala\par}
  {\itshape Associate Professor, 
  Department of Civil \& Environmental Engineering, \\
  University of Houston, Houston, Texas 77204. \\
  \textbf{phone:} +1-713-743-4418, \textbf{e-mail:} knakshatrala@uh.edu \\
  \textbf{website:} http://www.cive.uh.edu/faculty/nakshatrala}\\[0.1\baselineskip]

\begin{figure}[h]
  \centering
    \includegraphics[scale=0.24]{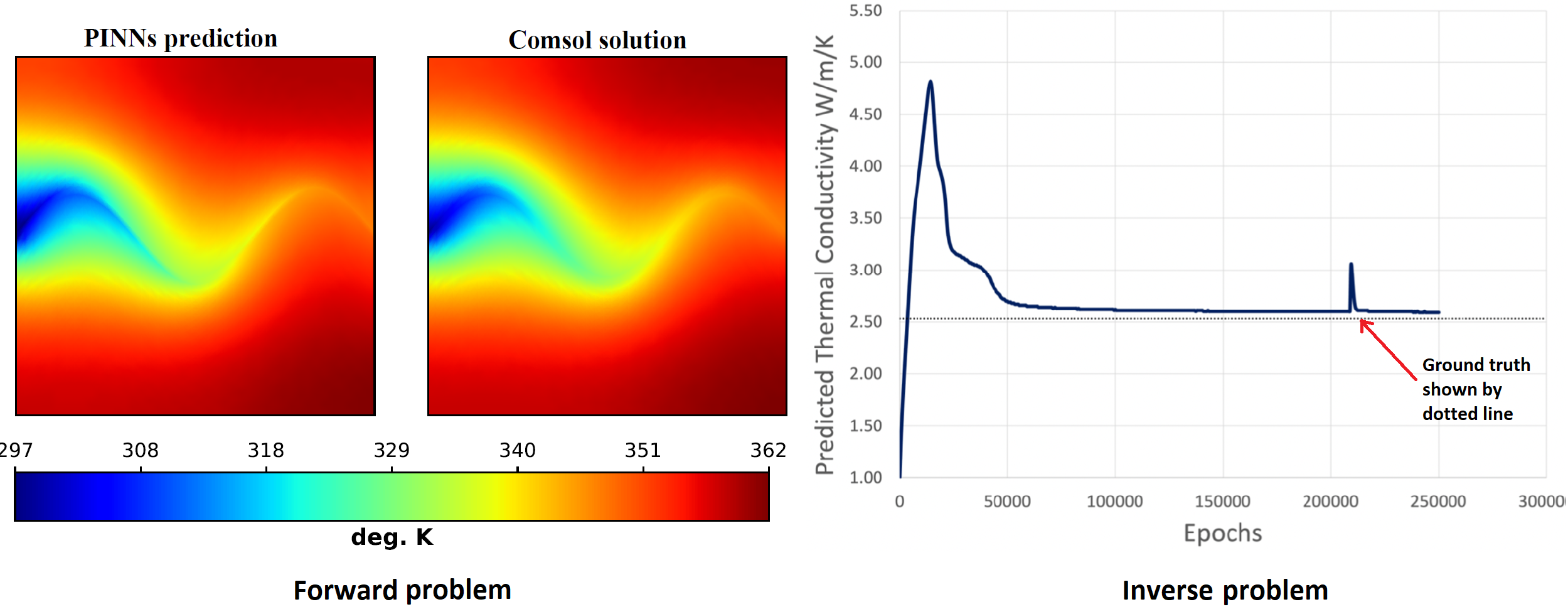}
    \emph{This figure shows the results from the CoolPINNs framework for forward and inverse problems. Heat is supplied to the bottom of the plate, and the top surface is free to convect and radiation. A flowing fluid through an embedded vasculature---in the shape of a sine wave---regulates the plate's temperature. The forward problem calculates the temperature field, whereas the inverse problem predicts thermal conductivity using noisy temperature data.}
\end{figure}
  
  \vfill
  {\scshape 2023} \\
  {\small Computational \& Applied Mechanics Laboratory} \par
\end{titlepage}

\begin{abstract}
Emerging technologies like hypersonic aircraft, space exploration vehicles, and batteries avail fluid circulation in embedded microvasculatures for efficient thermal regulation. Modeling is vital during the design and operational phases of these engineered systems. However, many challenges exist in developing a modeling framework. What is lacking is an accurate framework that (i) captures sharp jumps in the thermal flux across complex vasculature layouts, (ii) deals with oblique derivatives (involving tangential and normal components), (iii) handles nonlinearity because of radiative heat transfer, (iv) provides a high-speed forecast for real-time monitoring, and (v) facilitates robust inverse modeling. This paper addresses these challenges by availing the power of physics-informed neural networks (PINNs). We develop a fast, reliable, and accurate Scientific Machine Learning (SciML) framework for vascular-based thermal regulation---called CoolPINNs: a PINNs-based modeling framework for active cooling. The proposed mesh-less framework elegantly overcomes all the mentioned challenges. The significance of the reported research is multi-fold. \emph{First}, the framework is valuable for real-time monitoring of thermal regulatory systems because of rapid forecasting. \emph{Second}, researchers can address complex thermoregulation designs since the approach is meshless. \emph{Finally}, the framework facilitates systematic parameter identification and inverse modeling studies, perhaps the most significant utility of the current framework.
\end{abstract}

\maketitle

\vspace{-0.3in}

\section*{ABBREVIATIONS} 
\begin{tabular}{lr|lr} \hline
FEM & Finite Element Method 
& PDE & Partial Differential Equation \\ 
PINN & Physics-Informed Neural Network 
& ROM & Reduced-Order Model \\
SciML & Scientific Machine Learning 
& VTR & Vascular-based Thermal Regulation \\ \hline
\end{tabular}

\section{INTRODUCTION AND MOTIVATION}
\label{Sec:S1_PINNS_Intro}
Thermal regulation systems are used in diverse engineering applications like satellites \citep{lyall2008experimental},  space re-entry vehicles \citep{finke1990calculation}, photovoltaics \citep{brinkworth1997thermal}, geothermal heating and cooling systems \citep{rosen2017geothermal}, head cooling systems \citep{WOS:000646196300013}, 
bio-printing \citep{elomaa2017additive}, and extravehicular activity space suits \citep{nyberg2001model}. In all these cases, maintaining the temperature within a specific range is vital for the system's optimal performance.  For instance, the efficiency of the photovoltaic system reduces with an increase in system's mean surface temperature \citep{aldossary2016technical}. Likewise, low temperatures reduce capacity in lithium-ion batteries, while high temperatures degrade battery life \citep{PANCHAL2016204}.

An efficient method for thermal regulation is circulating a fluid through vasculatures embedded in the structural component. We refer to this method as \emph{vascular-based thermal regulation} (VTR)---the paper's primary focus. Prevalent in biological systems (e.g., cardiovascular system), VTR can maintain body temperature within the desired range, even when the surrounding temperature is considerably different \citep{hill1976jackrabbit, hall2020guyton}. Because of its ability to transport large amounts of heat compared to traditional methods, VTR is sought-after even in the synthetic world. 
Further, new fabrication techniques (e.g., deflagration of sacrificial components \citep{9595348} and vaporization of sacrificial components \citep{patrick2017robust}) have simplified the manufacturing of vasculatures in synthetic composite materials. Diverse applications of VTR include spacecrafts \citep{tan2018computational, ghosh2016small}, hypersonic vehicles \citep{gou2019design}, battery packs \citep{WOS:000413131200050}, multifunctional metamaterials \citep{devi2021microvascular}, photovoltaics \citep{huang2004thermal}, biomedical \citep{doi:10.1126/science.abl8532}, and shape memory alloys \citep{WOS:000419004600017}.

Many VTR applications require the vasculature cross-sections to be smaller than the component size to minimize their weight and volume.
Secondly, unless the vasculature volume fraction is very low, including a vascular network within a component substantially affects its mechanical properties like strength, stiffness, and interlaminar shear strength \citep{ phillips2011mechanical, WOS:000332440000002}. This is very critical for VTR applications like the casing of battery packs in electric vehicles, in which the composite panels containing vasculatures also provide structural support. Therefore, the hydraulic diameters of the vasculatures are often tiny compared to the panel dimensions \citep{bunce2017cubesat, phillips2011mechanical, kozola2010characterization}.
For these practical considerations, we focus on vasculatures with small hydraulic diameters in thin structural panels.
 
Modeling is irreplaceable during the design and operational phases of VTR applications.
Accurate modeling of VTR will require a solution of a 3D heat transfer equation coupled with the heat transfer by fluid flow inside the vasculature.
However, two simplifications in the underlying governing equations can be made for vasculatures with small hydraulic diameters in thin panels.
First, the temperature variation along thickness in thin panels will be small and a two-dimensional heat transfer model is sufficient. 
Second, the thermal load over the vasculature could be represented as a line load for vasculatures with a hydraulic diameter much smaller than the panel dimensions. 
With these two simplifications, the heat transferred by the fluid flow can be modeled by an \say{oblique derivative} boundary condition without explicitly modeling the fluid flow in the vasculature \citep{ramkrishna1979boundary}. The oblique derivative is untypical as it makes the boundary value problem non-self-adjoint \citep{ramkrishna1979boundary,ramkrishna1981boundary}.

The mathematical model resulting from these two simplifications closely approximates the coupled thermo-fluids governing equations and is referred to as a reduced order model (ROM)  \citep{nakshatrala2022qualitative}. 
The ROM is attractive since it reduces the problem dimensionality and eliminates the degrees of freedom associated with fluid flow.
However, approximation of the thermal load over a line leads to a discontinuity in thermal flux across the vasculature, modeled as a \say{jump term.} Therefore, the first requirement of the modeling framework is that the framework should be able to solve the equation with a jump term.

Another key requirement for the modeling framework is the ability to solve inverse problems for VTR applications with critical utility, like identifying material degradation of the vasculature panel by periodically calculating thermal conductivity using panel temperature data. However, solutions to inverse problems in heat transfer are generally ill-posed because their solution may become unstable due to small changes in the input parameters. Additionally, solving inverse problems requires time-consuming iterative solutions using the traditional methods \citep{ozisik2018inverse}. Therefore, the modeling framework to solve inverse problems for VTR applications should be robust and computationally efficient.

Modeling radiation heat transfer is essential for applications like equipment in space missions.
The inclusion of radiation using the Stefan–Boltzmann law makes the problem non-linear, and the framework should be able to solve non-linear equations. Optimization studies to achieve efficient heat transfer for critical operations lead to complex vasculature designs \citep{mcelroy2015optimisation}, which require very fine meshes using traditional numerical methods like FEM to provide accurate numerical solutions. To avoid meshing-related issues, a meshless method will be highly desirable.

Based on the above discussion, the modeling technique should (1) handle sharp jumps in the heat flux, (2) deal with inverse problems in an efficient and robust manner, (3) easily treat the non-linearity due to radiation, (4) be a mesh-less method for ease of modeling complex vasculatures, and (5) be fast and accurate. Currently, a modeling framework with the said features needs to be improved. However, the recent developments in the Scientific Machine Learning (SciML) methods have shown promise in developing a framework that can embody all these features.

A SciML method trains an algorithm to identify a pattern within extensive scientific data provided to it. Once trained, the algorithm can rapidly predict quantities of interest, such as solution fields. Over the past decade, SciML methods have grown exponentially and have been applied to a gamut of science and engineering problems. The accompanying reasons for this proliferation include the need to quickly process the massive data generated by sensors and simulations \citep{liu2019computational, reichstein2019deep}, increasing need to solve realistic simulations that are multiscale, involve multi-physics equations, and often need to be solved over large spatial and temporal domains \citep{OSS_SciML}, the ability of the trained ML models to make quick predictions, and advances in high-performance computing. Specific complex problems for which researchers have used SciML methods are fluid turbulence \citep{pandey2020perspective}, reactive-transport systems \citep{jagtap2022deep}, and the design of alloys \citep{durodola2022machine}, to name a few. 

One of the widely-used ML architectures at the forefront of developments in SciML is Deep Learning (DL) \citep{baker2019workshop}, a subset of ML that employs multiple hidden layers between a neural network's input and output. In a seminal paper by \citet{raissi2019physics}, DL is employed to develop \emph{Physics-informed Neural Networks} (PINNs). 
PINNs exploit the universal function approximation capability of deep neural networks and the automatic differentiation technique of scientific computing to constitute a loss term for the governing PDE and boundary conditions of the problem. 
This loss term imposes an additional constraint that limits the possible solutions. Therefore, the PINNs frameworks are efficient compared to the pure data science techniques of finding patterns in a vast amount of data. Furthermore, DL's automatic differentiation capability calculates the PDE's derivatives more efficiently and accurately than the traditional numerical differentiation \citep{baydin2018automatic}. Due to these advantages, the PINNs framework has been applied to solve various scientific problems since its publication in 2019 \citep{cuomo2022scientific, cai2022physics}. An overview of PINNs applied to heat transfer problems, which is the focus of the current paper, is presented below.

\citet{WOS:000680137900006} applied PINNs to solve the forward-posed problem for steady-state multi-species convective and diffusive flow and heat transfer problems. 
\citet{WOS:000646861400001} used PINNs for heat transfer problems with forced and mixed convection in the thermal design of power electronics. 
In \citep{WOS:000672478100022}, a PINNs methodology was developed to solve conduction heat transfer problems with convective heat transfer boundary conditions, which has applications in additive manufacturing and for parts heated in ovens. 
\citet{WOS:000700499500003} proposed PINNs for a reduced-order model without requiring the extra high-fidelity snapshots and demonstrated application to steady and unsteady natural convection problems. 
A PINNs framework for solving fluid problems in the rarefied and transitional regimes with applications in multiscale heat transfer was developed by \citet{WOS:000727776100012}. 
A PINNs framework to achieve accurate urban land surface temperature predictions was developed by \citet{WOS:000741407800002}.

All the aforementioned studies solved the PDEs without a \say{jump term} (discussed in \S\ref{Sec:S2_PINNS_GE}) or discontinuity in the underlying field.  On the other hand, very few attempts have been made to solve problems involving discontinuity in the underlying field using PINNs. 
\citet{WOS:000723617600011} developed a \say{control volume PINNs} framework for solving a Burger’s equation for a Riemann problem with zero viscosity. 
They employed regularization and artificial viscosity effects to model the rarefaction shock solution. \citet{WOS:000506874400038} have solved 1- and 2-dimensional Euler equations with a shock discontinuity for forward and inverse problems.

\begin{figure}
    \centering
    \includegraphics[scale=0.6]{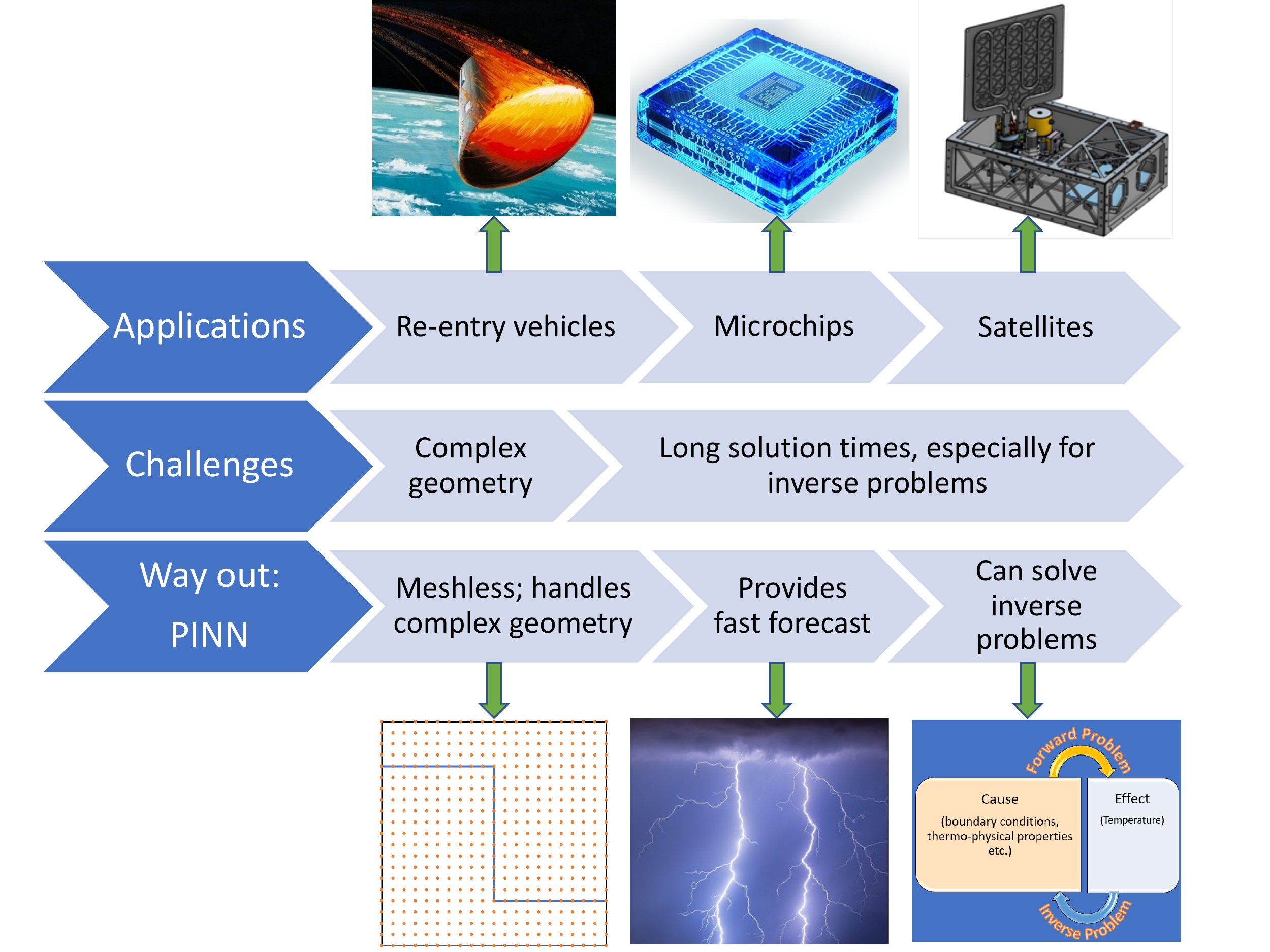}
    \caption{This figures outlines the applications and challenges of VTR and the advantages offered by modeling using PINNs. \label{fig:challenges_TR}}
\end{figure}

We present a PINNs framework meeting all the aforementioned requirements to solve a ROM for the VTR applications. 
We refer to this framework as CoolPINNs since we solve the problems for cooling applications, but the framework is directly applicable to heating applications.
We demonstrate that the CoolPINNs framework solves forward and inverse problems with underlying non-linear heat transfer equations with a jump term in a fast, robust, and accurate manner. Since our framework is PINNs-based, it is a meshless method and can easily model complex vasculature geometries. Forward problems for VTR applications have been studied using Finite Element Methods; e.g., see \citep{tan2016gradient, safdari2015nurbs, aragon2010generalized, nakshatrala2022qualitative, pejman2019gradient}. 
However, to the best of our knowledge, this is the first application of PINNs to solve forward and inverse heat transfer problems having a discontinuity in the thermal flux field. 
The challenges and benefits offered by PINNs for VTR applications are pictorially depicted in Fig.~\ref{fig:challenges_TR}.

The plan for the rest of this article is as follows. 
We first outline a reduced-order model for active cooling (\S\ref{Sec:S2_PINNS_GE}). We then present the architecture of CoolPINNs: the proposed PINNs-based modeling framework for active cooling (\S\ref{Sec:S3_PINNS_Framework}). Next, we illustrate, using representative results, the salient features of the modeling framework when solving forward problems (\S\ref{Sec:S4_PINNS_NR}). After that, we demonstrate the applicability of the proposed framework to inverse problems (\S\ref{Sec:S5_PINNS_INV}). Finally, the article concludes with concluding remarks alongside possible future works (\S\ref{Sec:S6_PINNS_Closure}).

\section{A MATHEMATICAL MODEL FOR VASCULAR-BASED THERMAL REGULATION}
\label{Sec:S2_PINNS_GE}

\subsection{Theory and background for the PDEs involving a jump term:}
Many partial differential equations (PDEs) describing common physical phenomena like fluid flow, material deformation, and heat transfer have the underlying kinematics and balance laws that assume a degree of smoothness in the fields such as density, strain, and heat flux. 
However, for certain real-life applications, the classical degree of smoothness is not achieved at a finite number of surfaces (for 3D) or lines (for 2D) within the domain. 
Such a surface or line is commonly called a \say{singular surface} or a \say{singular line.} In such cases, the underlying fields experience finite jump discontinuities at the singular surfaces but remain smooth in the rest of the domain.  
A typical example of a singular surface is the interface of the two dissimilar materials bonded together. 
Such a singular surface will have discontinuous stress and temperature gradients for structural and thermal loads. 
Other examples of singular surfaces include impact loading of solids, transonic flows in gas dynamics, and phase transitions in solids where the singular surface corresponds to a wavefront, a shock wave, and a phase boundary, respectively \citep{chadwick2012continuum,abeyaratne1998continuum}.

For the problem at hand, modeling of the heat carried by the fluid in the vasculature as a thermal line load leads to a discontinuity or a \emph{jump term} in the thermal flux across the vasculature. 
The jump term includes an \emph{oblique derivative} \citep{ramkrishna1981boundary, WOS:000251903300002}---in terms of the temperature field. Meaning, the direction of the gradient makes an oblique angle to the direction of the vasculature. Thus, the jump term will combine normal derivative and tangential derivative components. 

\subsection{Problem description and governing equations}

\begin{figure}
  \centering
    {\includegraphics[scale=0.4]{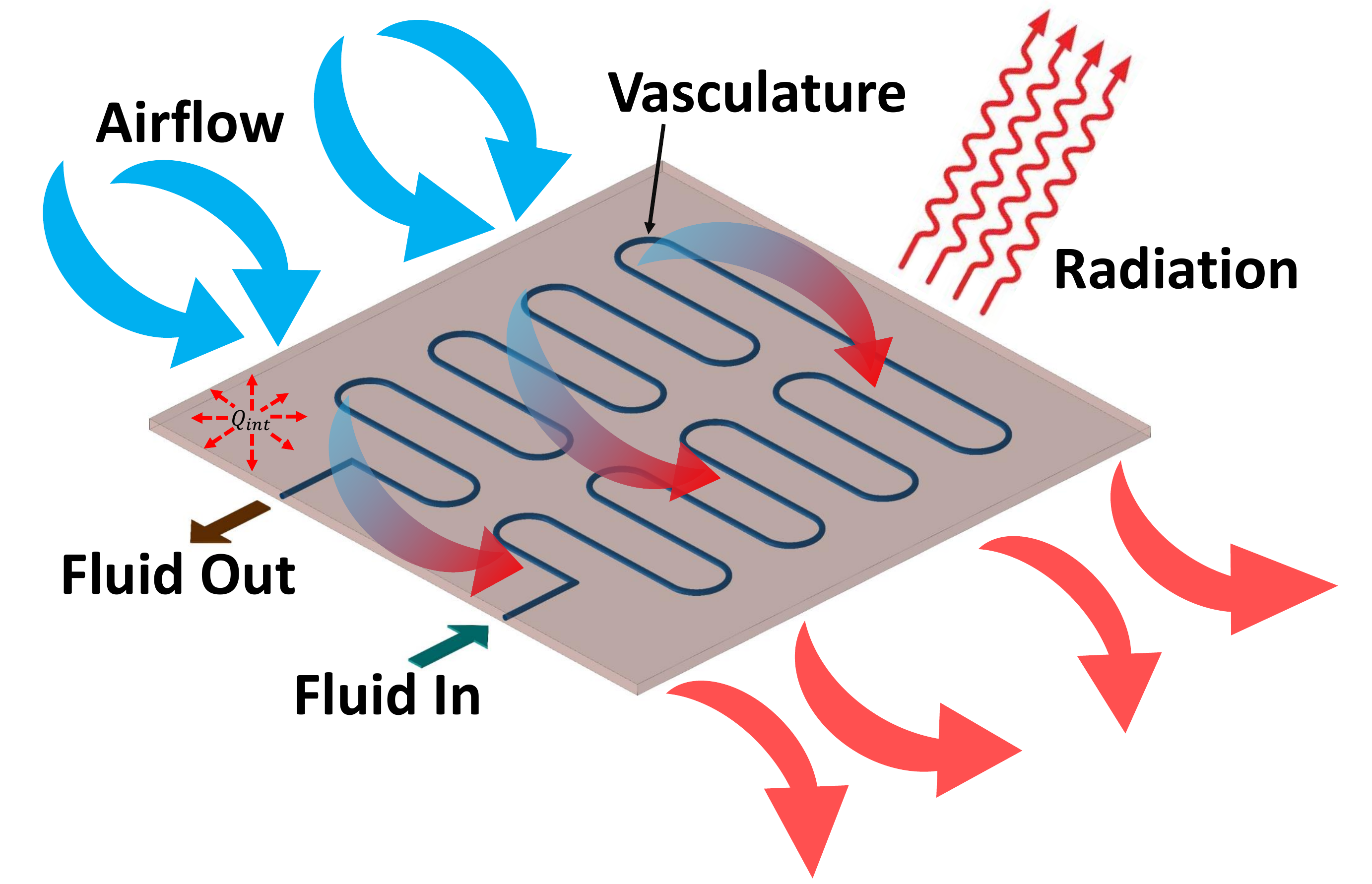}}
    \caption{The figure illustrates a typical active-cooling application. The system comprises a thin host solid---idealized as a plate---with an embedded vasculature. An external source supplied heat to the system. Within the body, the transfer of heat occurs via conduction. At the same time, heat is transferred to surroundings by convection and radiation. A fluid flows through the vasculature to regulate the plate's temperature.}
    \label{fig:AC_Illustration}  
\end{figure}

Consider a thin flat plate as shown in Fig. \ref{fig:AC_Illustration} with a thickness \emph{t} much smaller than its lateral dimensions. The plate has an internal source that generates heat energy. 
In practical applications, such a plate could be attached to a component generating heat, and as such, heat may not be generated within the plate itself. 
Heat is conducted within the plate and is transferred to the surroundings through convection and radiation from the top surface. 
A vasculature is embedded at the mid-surface of the plate and has a tiny hydraulic diameter compared to the lateral dimensions of the plate. 
The vasculature traverses through the plate, separating the plate mid-surface into two separate sub-domains. 
A fluid is circulated through the vasculature to regulate the plate temperature. 

To derive the governing equations for the thermal regulation with active cooling, we start with steady-state heat conduction in a three-dimensional body $\mathcal{B}$ with a heat-generating source. We assume the Fourier model, which takes the following form: 
\begin{align}
    \mathbf{q}(\mathbf{x}) = - \mathrm{K}(\mathbf{x}) \, \mathrm{grad}[\theta] 
    \quad \mathrm{in} \; \mathcal{B}, 
\end{align}
where $\mathbf{x}$ is the position vector, $\mathrm{\theta}(\mathbf{x})$ is the temperature field, $\mathbf{q}(\mathbf{x})$ the heat flux vector field, $\mathrm{K}(\mathbf{x})$ the thermal conductivity, and $\mathrm{grad}[\bullet]$ is the spatial gradient operator. The balance of energy reads: 
\begin{align}
    \label{Eqn:SSHeat}
    -\mathrm{div}[\mathbf{q}(\mathbf{x})] + F(\mathbf{x}) = 0 \;\;\;\; \text{or} \;\;\;\;
    \mathrm{div}[\mathrm{K(\mathbf{x})}\,\mathrm{grad}[\theta(\mathbf{x})]] + F(\mathbf{x}) = 0 \quad \mathrm{in} \; \mathcal{B},
\end{align}
where div[$\bullet$] represents the divergence operator, and F($\mathbf{x}$) is the rate of heat generation per unit volume. 

Since the plate is thin, it is reasonable to assume that the temperature along the plate thickness direction is constant, and hence a two-dimensional reduced-order model is sufficient. To this end, we write the governing equation on the mid-plane of the body---which we refer to as the domain $\Omega$.  Assuming that the plate material has a constant thermal conductivity and heat is generated uniformly over the plate area, Eq.~\eqref{Eqn:SSHeat} can be simplified as:
\begin{align}
    \label{Eqn:SSHeat2D}
    t\,\mathrm{K}\,\mathrm{div}[\mathrm{grad}[\theta(\mathbf{x})]] + f(\mathbf{x}) = 0 \;\;\;\;\text{or}\;\;\;\; t\,\mathrm{K}\,\nabla^{2}\theta(\mathbf{x}) + f(\mathbf{x}) = 0 \quad \mathrm{in} \; \Omega,
\end{align}
where $t$ is the plate thickness, $f(\mathbf{x}) \,(=t \, F(\mathbf{x}))$ is the steady-state rate of heat generation per unit area, and $\nabla^{2}$ is the Laplacian operator. In the above reduced-order model, $\mathbf{x}$ is the 2D position vector.

To account for the heat transferred to the surroundings through convection and radiation, we modify Eq.~\eqref{Eqn:SSHeat2D} as follows:  
\begin{align}
    \label{Eqn:SSHt_K_h_rad}
    t\,\mathrm{K}\,\nabla^{2}\theta(\mathbf{x}) + f(\mathbf{x}) - h_{T} \, (\theta(\mathbf{x}) - \theta_{\mathrm{amb}}) - \varepsilon \, \sigma \,  (\theta^{4}(\mathbf{x}) - \theta_{\mathrm{amb}}^{4}) = 0 \quad \mathrm{in} \; \Omega,
\end{align}
where $h_T$ is the convective heat transfer coefficient,  $\theta_{\mathrm{amb}}$ is the ambient temperature, $\varepsilon$ is the emissivity of the plate material, $\sigma$ is the Stefan-Boltzmann constant. Equation~\eqref{Eqn:SSHt_K_h_rad} needs to be augmented by a jump condition along the vasculature to account for the heat transported by the flowing fluid in the vasculature.

\begin{figure}
    \centering
        {\includegraphics[scale=0.4]{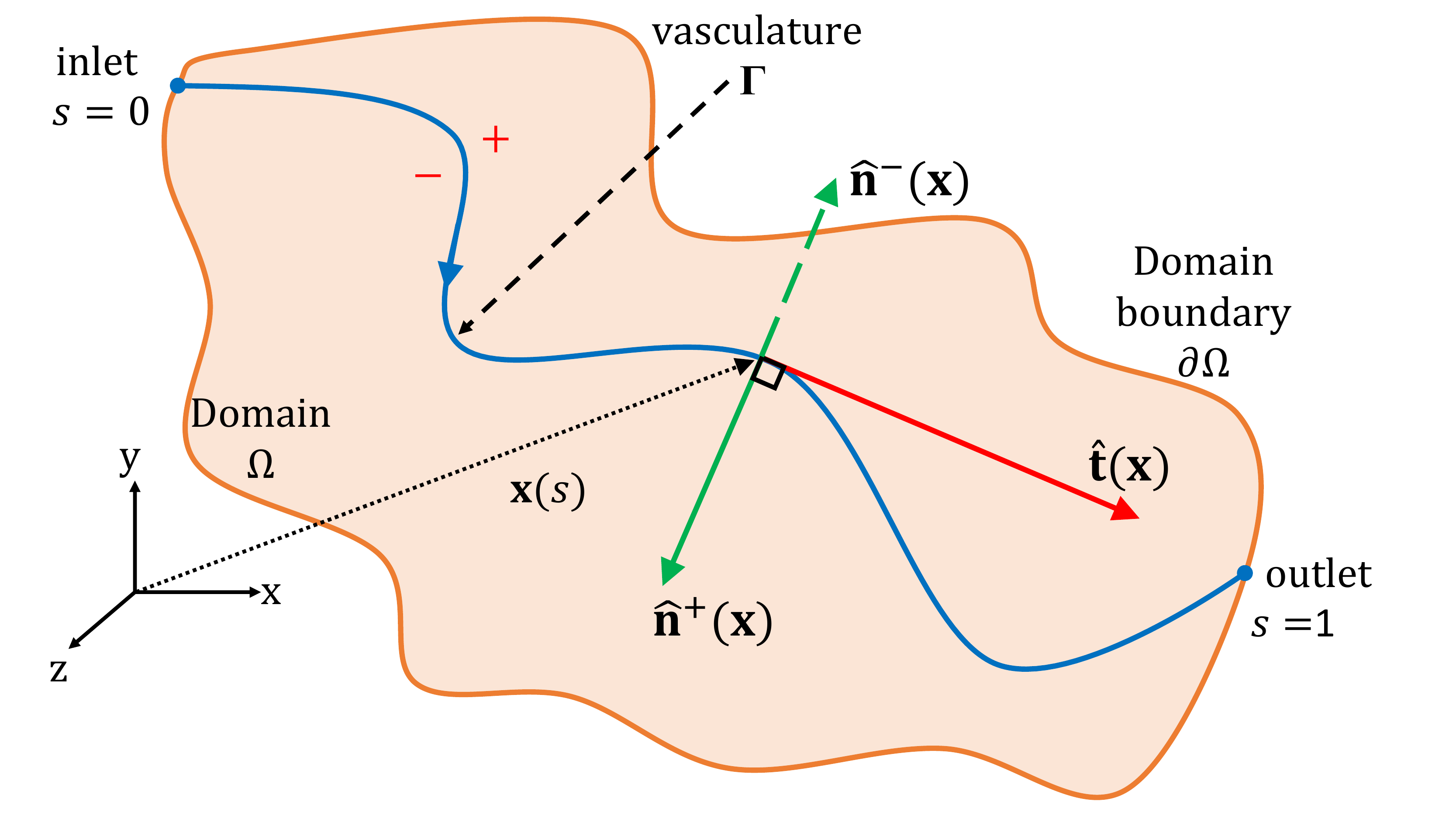}}
        \caption{An arbitrary-shaped domain, a vasculature, and the nomenclature used in the equations.}
        \label{fig:vasc_schematic}  
\end{figure}

Figure \ref{fig:vasc_schematic}
shows an arbitrary vasculature, denoted as $\Gamma$. 
Let us denote the two sub-domains separated by the vasculature using $+$ and $-$ signs. The vasculature geometry is defined by a position vector $\mathbf{x}(s)$ with $s$ denoting the arc-length: $s = 0$ and $s = 1$ denote the inlet and outlet of $\Gamma$, respectively. 
The unit tangent vector to $\Gamma$ along the fluid flow direction is $\widehat{\mathbf{t}}(\mathbf{x})$. 
The unit normal vector to $\Gamma$ for the sub-domains $+$ and $-$ are $\widehat{\mathbf{n}}^{+}(\mathbf{x})$ and $\widehat{\mathbf{n}}^{-}(\mathbf{x})$. 
Note that $\widehat{\mathbf{n}}^{+}(\mathbf{x})$ = --$\widehat{\mathbf{n}}^{-}(\mathbf{x})$.

With this information, we proceed to define the equations for the \say{jump term} in the thermal flux. 
Note that the temperature field is continuous across the vasculature, but the heat flux (or temperature gradient) is discontinuous across the vasculature. 
The jump term for heat flux at any point $\mathbf{x}$ on the vasculature is defined as
\begin{align}
    \label{Eqn:DefFluxJump}
    [\![ \mathbf{q}(\mathbf{x}) ]\!] = \mathbf{q}^{+}(\mathbf{x})- \mathbf{q}^{-}(\mathbf{x}),
\end{align}
where $\mathbf{q}^{+}(\mathbf{x})$ and $\mathbf{q}^{-}(\mathbf{x})$ are the limiting values of the heat flux at point $\mathbf{x}$ on the vasculature as we approach $\mathbf{x}$ from the $+$ and $-$ sub-domains along the direction perpendicular to the vasculature.

\begin{figure}
    \centering
    \includegraphics[scale=0.5]{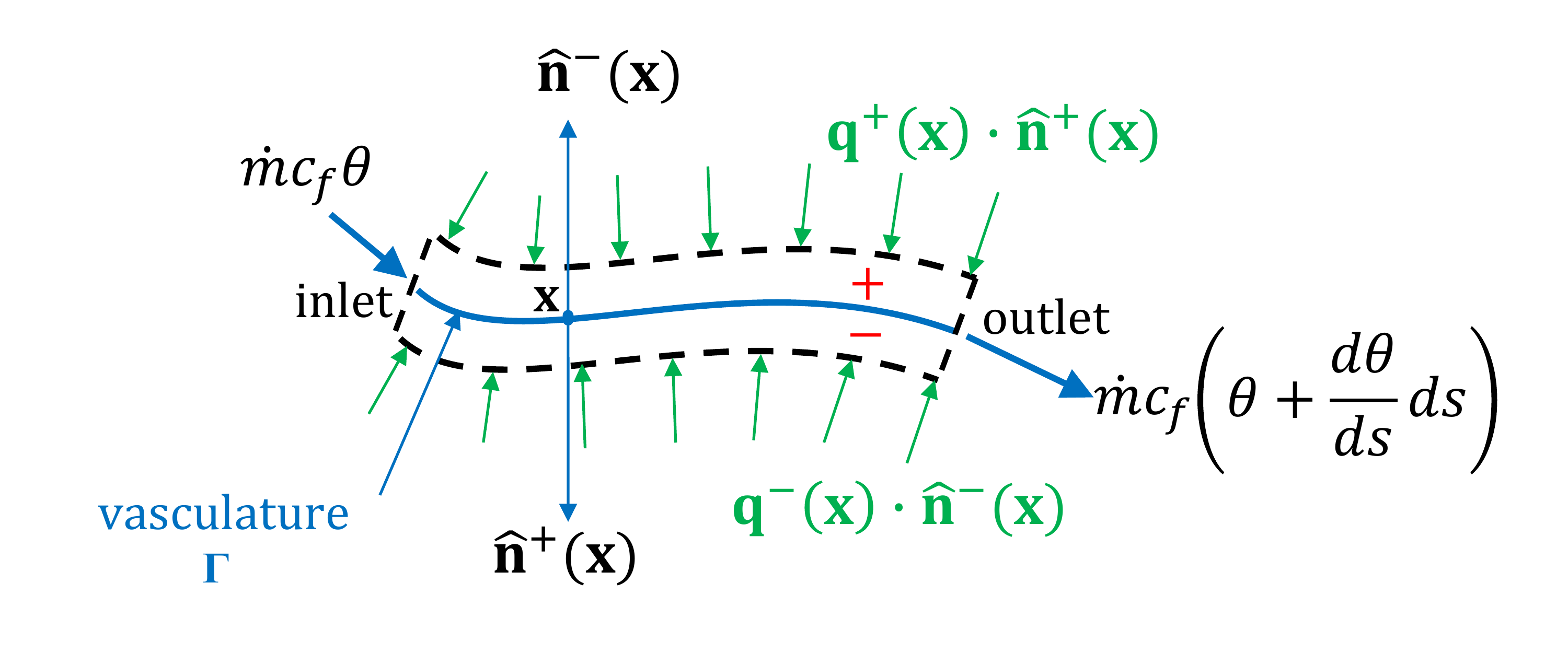}
    \caption{Heat energy balance along a small vasculature segment.}
    \label{fig:InfiVasc}
\end{figure}

To derive the mathematical formula for the jump term in heat flux, we consider an infinitesimal vasculature segment as shown in Fig. \ref{fig:InfiVasc}. 
Applying the law of conservation of energy to this segment, we get
\begin{align}
    \label{Eqn:BalHt}
    \dot{m} \, c_{f} \, \theta 
    + \Big(\mathbf{q}^{+}(\mathbf{x})\cdot\widehat{\mathbf{n}}^{+}(\mathbf{x}) + \mathbf{q}^{-}(\mathbf{x}\Big)\cdot\widehat{\mathbf{n}}^{-}(\mathbf{x}))\,t\,ds = \dot{m} \, c_{f} \, \left(\theta + \frac{d\theta}{ds}ds\right), 
\end{align}
where $\dot{m}$ is the mass flow rate of the fluid flow through the vasculature, $c_{f}$ is the specific heat capacity of the fluid, and $\dot{m}c_{f}\theta$ and $\dot{m}c_{f}(\theta + \frac{d\theta}{ds}ds)$ are the rate at which heat is carried in and out of the vasculature. 
Cancelling the common terms on the two sides of the equation and using the relation $\widehat{\mathbf{n}}^{+}(\mathbf{x})$ = $-\widehat{\mathbf{n}}^{-}(\mathbf{x})$, we get
\begin{align}
    \label{Eqn:BalHt2}
    t\,\Big(\mathbf{q}^{+}(\mathbf{x}) - \mathbf{q}^{-}(\mathbf{x})\Big)\cdot\widehat{\mathbf{n}}^{+}(\mathbf{x}) = \dot{m}c_{f}\frac{d\theta}{ds}.
\end{align}

Substituting the expression for the heat flux jump term from Eq.~\eqref{Eqn:DefFluxJump}, we get the following expression
\begin{align}
    \label{Eqn:FluxJumpFormula}
    t\,[\![ \mathbf{q} ]\!]\cdot\widehat{\mathbf{n}}^{+}(\mathbf{x}) = \dot{m} \, c_{f} \, \frac{d\theta}{ds} \;\;\;\;\;\; \text{or} \;\;\;\;\;\; t\,[\![ \mathbf{q} ]\!]\cdot\widehat{\mathbf{n}}^{+}(\mathbf{x}) = \dot{m} \, c_{f}\, \mathrm{grad}[\theta]\cdot\widehat{\mathbf{t}}(\mathbf{x}).
\end{align}

Equation \eqref{Eqn:FluxJumpFormula} is referred to as oblique derivative boundary condition in the literature \citep{ramkrishna1981boundary, WOS:000251903300002} because the resultant vector of the two vector terms in Eq.~\eqref{Eqn:FluxJumpFormula} is along a direction oblique to the normal (or tangent) to the boundary. 

For defining the boundary value problem, we  denote the two-dimensional mid-surface passing through the vasculature by $\Omega$. 
The boundary of this domain will be referred to as $\partial\Omega$. The domain boundary $\partial\Omega$ is divided such that $\partial\Omega$ = $\partial\Omega^{\theta}$ $\cup$ $\partial\Omega^{q}$ and $\partial\Omega^{\theta}$ $\cap$ $\partial\Omega^{q}$ = $\emptyset$, where $\partial\Omega^{\theta}$ and $\partial\Omega^{q}$ are portions of the boundary over which Dirichlet (temperature) and Neumann (heat flux) boundary conditions are specified.

The above governing equations for the reduced-order model are summarized below: 
\begin{subequations}
\label{Eqn:GE}
\begin{alignat}{2}
    \label{Eqn:BoE} 
    &t\,\mathrm{K}\,\nabla^{2}\theta(\mathbf{x}) + f(\mathbf{x}) 
    - h_{T} \, (\theta(\mathbf{x}) - \theta_{\mathrm{amb}}) 
    - \varepsilon \, \sigma \,  (\theta^{4}(\mathbf{x}) - \theta_{\mathrm{amb}}^{4}) = 0
    && \quad \mathrm{in} \; \Omega, \\
    \label{Eqn:q_jump_condition} 
    &t\,[\![ \mathbf{q} ]\!]\cdot\widehat{\mathbf{n}}^{+}(\mathbf{x}) = \dot{m} \, c_{f}\, \mathrm{grad}[\theta]\cdot\widehat{\mathbf{t}}(\mathbf{x}) 
    && \quad \mathrm{on} \; \Gamma, \\
     \label{Eqn:q_BC} 
    & \mathbf{q}(\mathbf{x}) \cdot \widehat{\mathbf{n}}(\mathbf{x}) 
    = q^{\mathrm{p}}(\mathbf{x}) 
    && \quad \mathrm{on} \; \partial\Omega^{q},  \\
    \label{Eqn:temp_BC} 
    &\theta(\mathbf{x}) = \theta^{\mathrm{p}}(\mathbf{x}) 
    && \quad \mathrm{on} \; \partial\Omega^{\theta}, 
    \, \mathrm{and} \\
    \label{Eqn:temp_BC_V_in} 
    &\theta = \theta_{\mathrm{in}}
    && \quad \mathrm{at} \; s = 0. 
\end{alignat}
\end{subequations}
In above equations, $t$ is the plate thickness, K is the thermal conductivity of the plate material, $\theta(\mathbf{x})$ is the temperature field, $f(\mathbf{x})$ is the heat source term, $h_T$ is the convective heat transfer coefficient,  $\theta_{\mathrm{amb}}$ is the ambient temperature, $\varepsilon$ is the emissivity, $\sigma$ is the Stefan-Boltzmann constant, $[\![ \mathbf{q}(\mathbf{x}) ]\!]$ is the jump term for heat flux across $\Gamma$, $\widehat{\mathbf{n}}^{+}(\mathbf{x})$ is the unit normal vector to sub-domain $+$ at $\Gamma$, $\dot{m}$ is the mass flow rate of the fluid flow through the vasculature, $c_{f}$ is the specific heat capacity of the fluid,  $\widehat{\mathbf{t}}(\mathbf{x})$ is the unit tangent vector to $\Gamma$ along the fluid flow direction, $\mathbf{q}(\mathbf{x})$ is the heat flux, $\theta^{\mathrm{p}}(\mathbf{x})$ and $\mathbf{q}^{\mathrm{p}}(\mathbf{x})$ are the prescribed temperature and flux on the boundary, and $\theta_{\mathrm{in}}$ is the temperature of the coolant at the inlet.

Because of the radiation term, the mathematical model is nonlinear. Also, although the current discussion focuses on active cooling applications, this setup applies even to active heating applications. For more details on this reduced-order model, see \citep{nakshatrala2022qualitative}. 

\section{PROPOSED MODELING FRAMEWORK: CoolPINNs}
\label{Sec:S3_PINNS_Framework}

\begin{figure}
    \centering
    \includegraphics[scale=0.42]{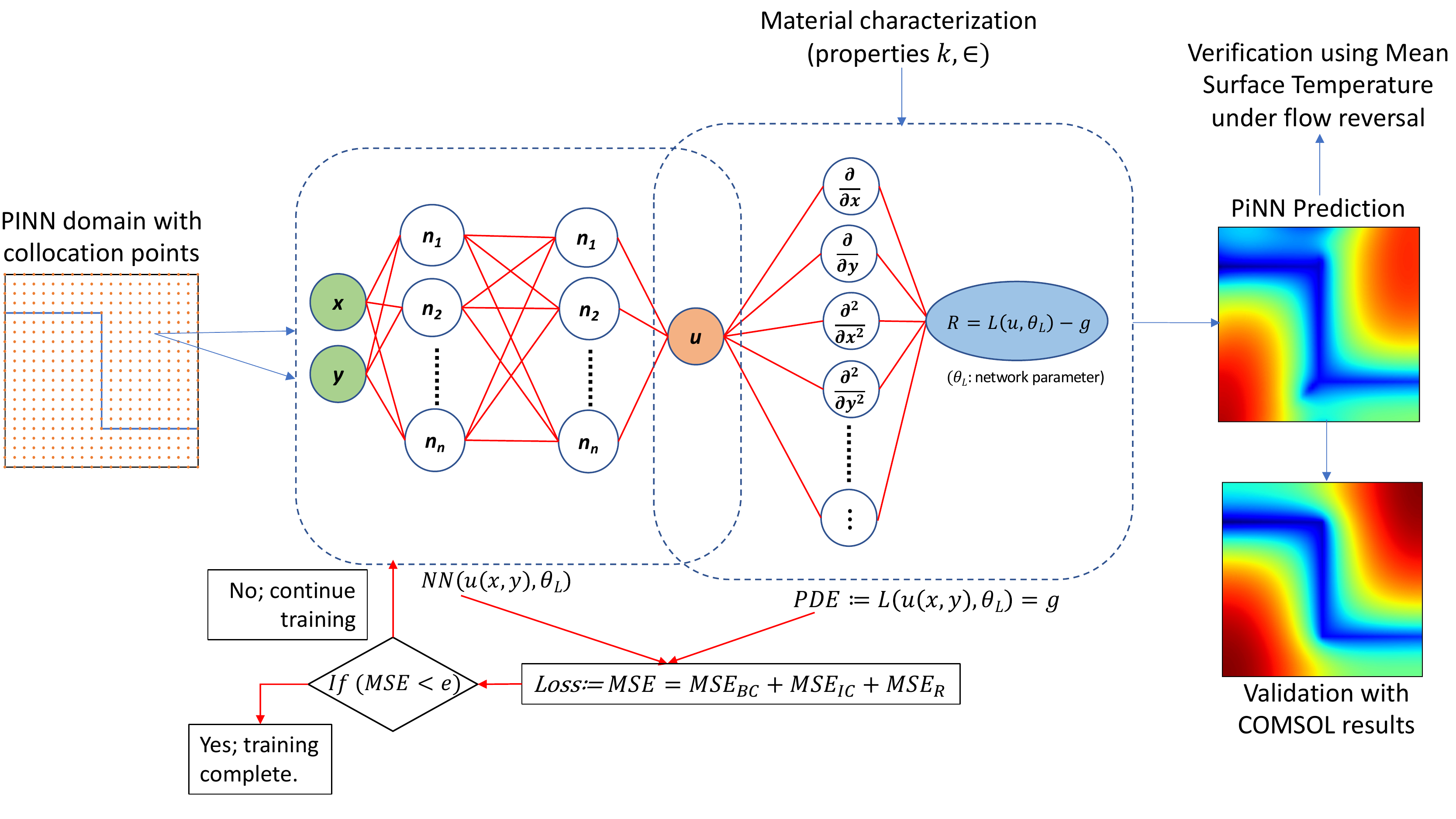}
    \caption{A schematic showing workflow of CoolPINNs framework to thermal regulation problems. 
    The CoolPINNs framework consists of two distinct neural networks shown by the dotted lines. 
    These two networks share the same hyperparameters and act together to minimize the loss function. 
    The network denoted as NN contributes to the loss function for the boundary training data. 
    In contrast, the network denoted as PDE contributes to the loss for the given thermal transport PDE and brings the physics-informed aspect to the architecture. 
    The resulting predictions from the framework are verified by comparing the results for mean surface temperature under flow reversal conditions and are validated by comparing with the FEM solutions obtained from the \textsf{Comsol} multi-physics package.
    \label{fig:PiNN_framework}}
\end{figure}

The PINNs framework developed herein---addressing vascular-based active-cooling applications---is called CoolPINNs. Nevertheless, the approach is applicable even for general thermal regulation applications, including heating. The deep learning architecture of the proposed framework comprises two separate neural networks, as shown by the dotted lines in Fig. \ref{fig:PiNN_framework}. These two networks share the same hyperparameters and act together to minimize a loss function, which needs to be constructed carefully. 
The network denoted as \textit{NN} contributes to the loss based on the boundary training data. On the other hand, the network represented by PDE contributes to the loss for the given PDE and brings the physics-informed aspect into the deep learning neural architecture. 
The PDE network utilizes the automatic differentiation of neural networks for the spatial dimensions and model parameters \citep{raissi2019physics,baydin2018automatic}. 
The two networks \textit{NN} and PDE combinedly ``learns'' and ``forecasts'' the temperature field $\theta_{\mathrm{L}}(\mathbf{x})$.

To define the loss function for the PDE network, we introduce $\mathrm{f}_{\mathrm{n}}(\mathbf{x})$ obtained by replacing $\theta(\mathbf{x})$ in Eq.~\eqref{Eqn:BoE} with $\theta_{\mathrm{L}}(\mathbf{x})$:
\begin{equation}
    \label{Eqn:BoE_RA} 
    \mathrm{f_{n}(\mathbf{x})\,:= t\,\mathrm{K}\,\nabla^{2}\theta_{\mathrm{L}}(\mathbf{x})} + f(\mathbf{x}) 
    - h_{T} (\theta_{\mathrm{L}}(\mathbf{x}) - \theta_{\mathrm{amb}}) 
    - \varepsilon \, \sigma \,  (\theta_{\mathrm{L}}^{4}(\mathbf{x}) - \theta_{\mathrm{amb}}^{4}).
\end{equation}
The reason for introducing a new variable is to distinguish between the solution from the neural network ($\theta_{\mathrm{L}}(\mathbf{x})$) \emph{versus} the exact solution of the governing equations \eqref{Eqn:BoE}--\eqref{Eqn:temp_BC_V_in}, denoted by $\theta(\mathbf{x})$. 

\subsection{Forward problems}
For forward problems, the governing PDE and boundary conditions are known \emph{a priori}; both are part of the description of the boundary value problem. So, the loss function consists of two parts. The loss contribution from the PDE (Loss\textsubscript{PDE}) represents the error residual at the collocation points within the domain and is defined as follows:
\begin{equation}
    \label{Eqn:PDELoss}
    \mathrm{Loss_{PDE} = \frac{1}{N_{PDE}}\sum_{i=1}^{N_{PDE}} \left| f_{n}(\mathbf{x}_{PDE}^{i}) \right|^{2}},
\end{equation}
where N\textsubscript{PDE} are the number of collocation points defined within the domain, and $\mathbf{x}_{\mathrm{PDE}}^{\mathrm{i}}$ are the coordinates of the collocation points within the domain. Further, it is necessary that the training data is available to account for the boundary conditions: that is, the temperature field for the collocation points on the Dirichlet part of the boundary ($\partial \Omega^{\theta}$) and the normal component of the heat flux for the collocation points on the Neumann part of the boundary ($\partial \Omega^{q}$). 
The loss on the prescribed boundary conditions (Loss\textsubscript{BC}) is defined as
\begin{equation}
    \label{Eqn:BCLoss}
    \mathrm{Loss_{BC} = \frac{1}{N_{BCD}}\sum_{i=1}^{N_{BCD}} \left| \theta_{BC}(\mathbf{x}^{i}) - \theta_{L}(\mathbf{x}^{i})\right|^{2} + \frac{1}{N_{BCN}}\sum_{i=1}^{N_{BCN}} \left| q_{BC}(\mathbf{x}^{i}) - q_{L}(\mathbf{x}^{i})\right|^{2}},
\end{equation}
where $\mathrm{N_{BCD}}$ and $\mathrm{N_{BCN}}$ are the number of boundary points over which Dirichlet and Neumann boundary conditions are specified, $\mathbf{x}^{\mathrm{i}}$ denotes a collocation point on the boundary, $\mathrm{\theta_{BC}}$ and $\mathrm{q_{BC}}$ are the prescribed Dirichlet and Neumann boundary conditions at $\mathbf{x}^{\mathrm{i}}$, and $\mathrm{q_{L}}(\mathbf{x}^{\mathrm{i}})$ is the thermal flux at $\mathbf{x}^{\mathrm{i}}$. Mathematically, 
\begin{equation}
    q_{\mathrm{L}}(\mathbf{x}^{\mathrm{i}}) = - \mathrm{K} \, \mathrm{grad}[\theta_{\mathrm{L}}] \cdot \widehat{\mathbf{n}}(\mathbf{x}) \Big|_{\mathbf{x} = \mathbf{x}^{\mathrm{i}}}
\end{equation}
Note that the deep neural network solves for $\theta_{\mathrm{L}}(\mathbf{x})$ at all the collocation points, which are on the domain as well as on the Dirichlet and Neumann parts of the boundary. Thus, the total loss that is minimized for forward problems during the training process of PINNs---denoted as $\mathrm{Loss_{fw}}$---is the sum of the above two losses (given by Eqs.~\eqref{Eqn:PDELoss} and \eqref{Eqn:BCLoss}):
\begin{equation}
    \label{Eqn:Loss} 
    \mathrm{Loss_{fw} = Loss_{PDE} + Loss_{BC}}.
\end{equation}

\subsection{Inverse problems}
For inverse problems, the general form of governing equation and nature of boundary conditions are known except for a few parameters. These parameters, which need to be predicted, can be a material property (e.g., thermal conductivity, which is the case in this paper) or even boundary conditions. In order to solve an inverse problem, we need to know solution (i.e., temperature or heat flux) over a small set of points within the domain; in this paper, we assume that the temperature is known at a sub-collection of collocation points with the domain. 

The loss associated with the PDE and BCs are calculated using Eqs.~ \eqref{Eqn:PDELoss} and \eqref{Eqn:BCLoss}.
The loss associated with the temperature solution over a small set of points within the domain is defined as follows:
\begin{equation}
    \label{Eqn:SolnLoss}
    \mathrm{Loss_{U} = \frac{1}{N_{u}}\sum_{i=1}^{N_{u}} \left| \theta_{u}(\mathbf{x}^{i}) - \theta_{L}(\mathbf{x}^{i})\right|^{2}},
\end{equation}
where $\mathrm{N_{u}}$ is the number of points used for training where temperature solution is known, and $\theta_{\mathrm{u}}$, $\theta_{\mathrm{L}}$ are the known and predicted temperatures at the points $\mathbf{x}^{\mathrm{i}}$.
The total loss minimized for inverse problems during the training process of PINNs is

\begin{equation}
    \label{Eqn:Loss_inv} 
    \mathrm{Loss_{inverse} = Loss_{PDE} + Loss_{BC} + Loss_{U}}.
\end{equation}

\begin{figure}
  \centering  
  \begin{subfigure}[b]{0.4 \textwidth}
    {\includegraphics[width=\textwidth]{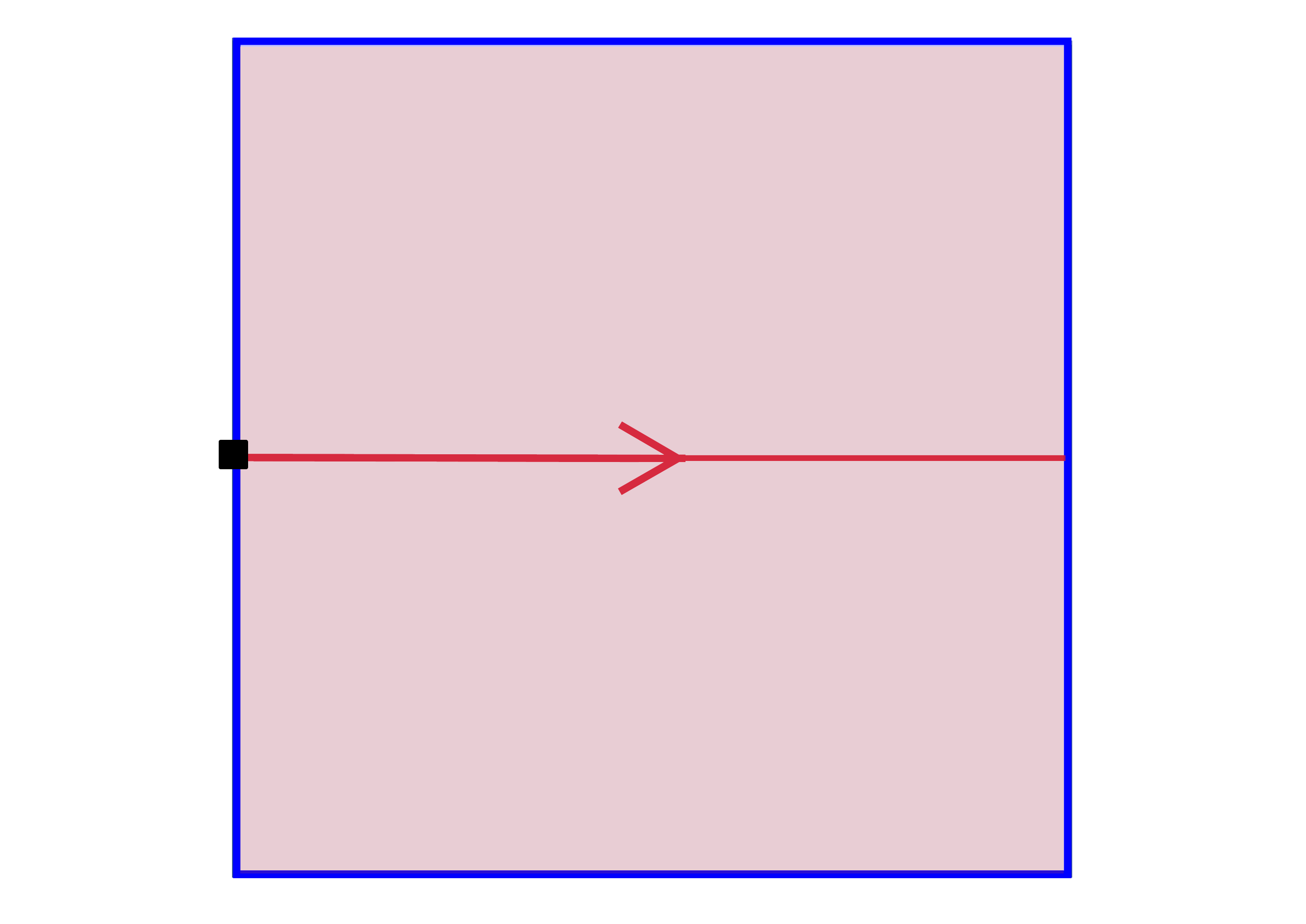}}
    \caption{\tiny Test problem \#1: Straight channel at center \\ of a square domain}
  \end{subfigure}
  %
  %
  \begin{subfigure}[b]{0.4 \textwidth}
    {\includegraphics[width=\textwidth]{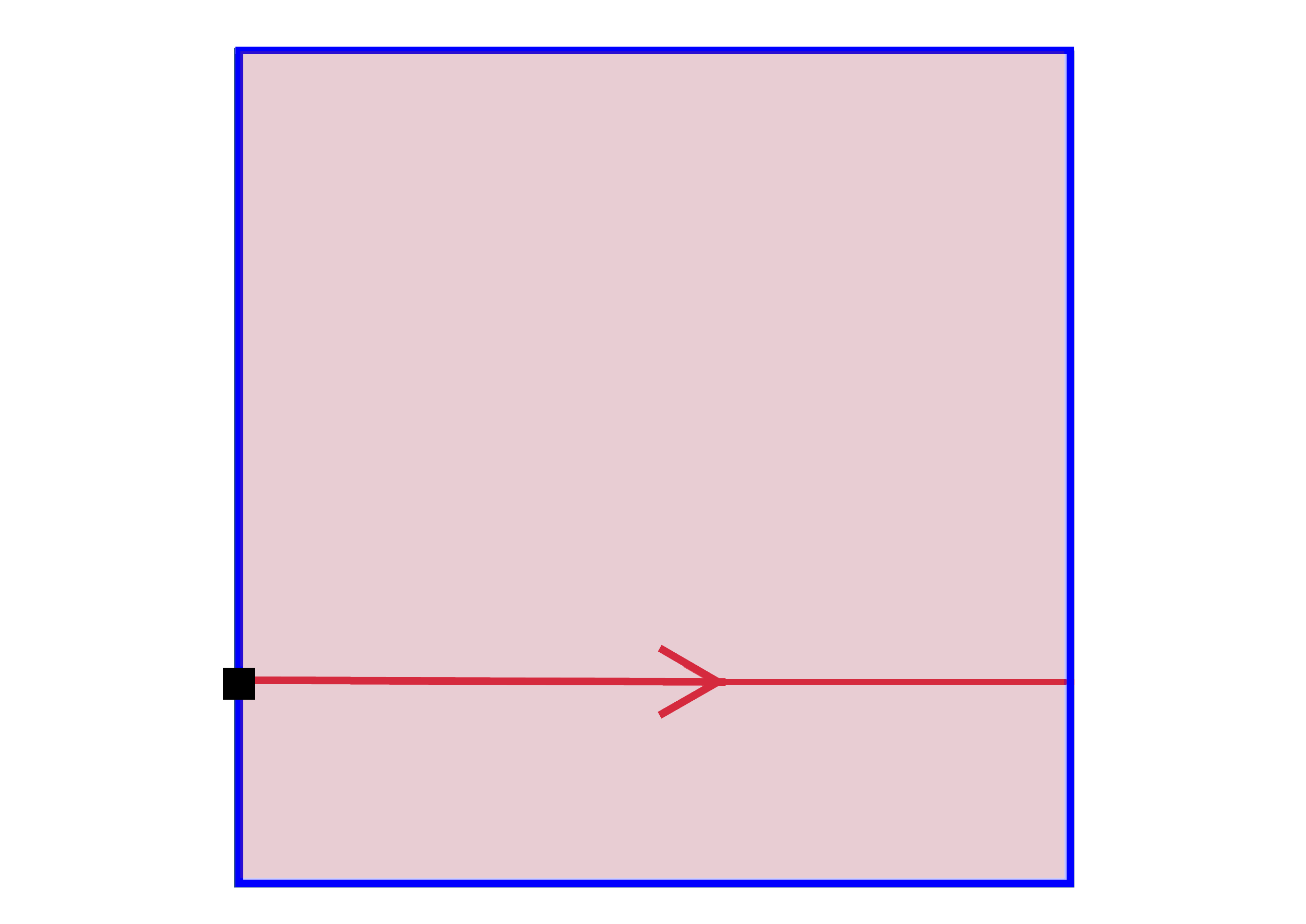}}
    \caption{\tiny Test problem \#2: Straight channel at off-center \\ of a square domain}
  \end{subfigure}
  %
  %
  \begin{subfigure}[b]{0.4 \textwidth}
    {\includegraphics[width=\textwidth]{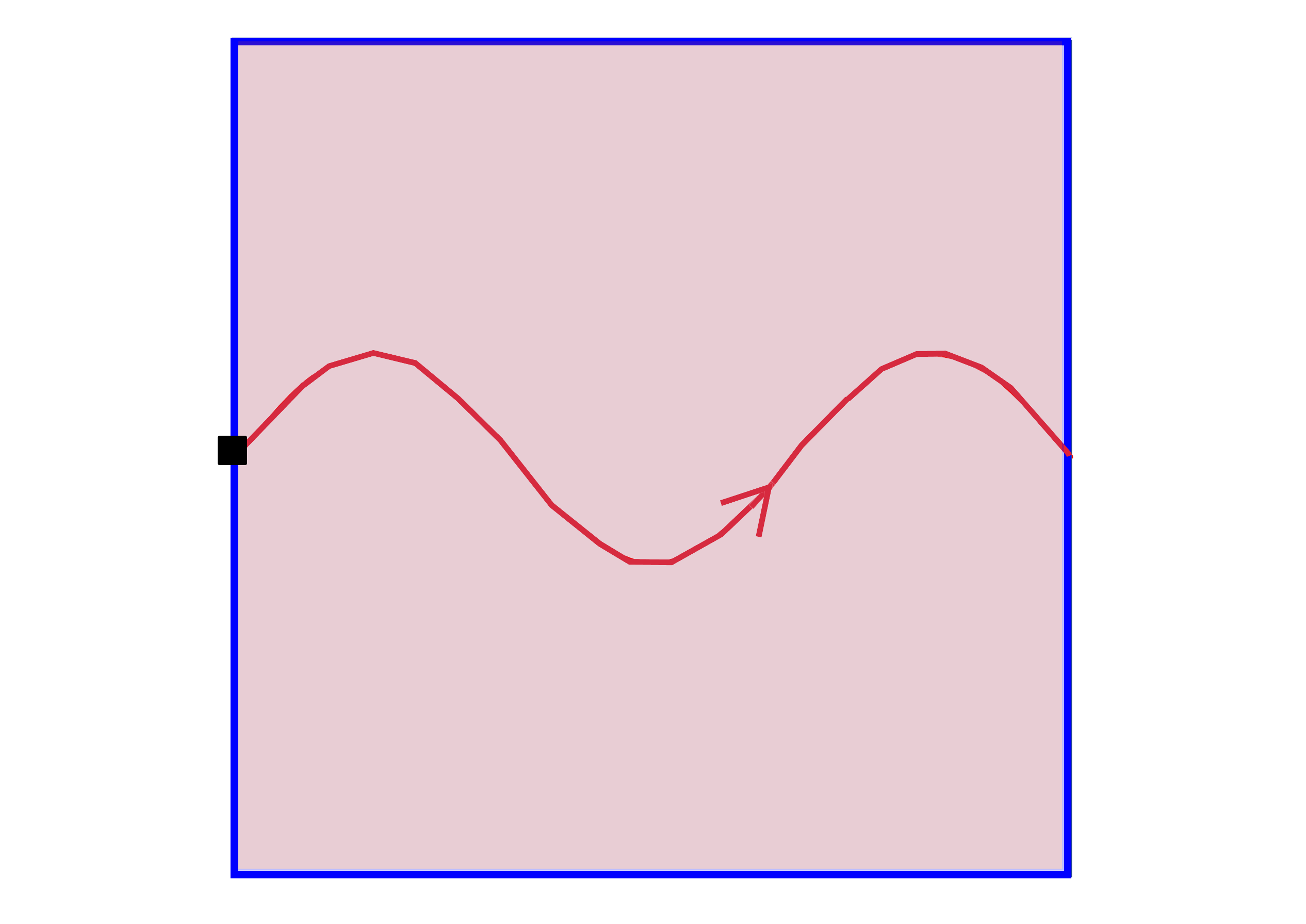}}
    \caption{\tiny Test problem \#3: Sine wave channel through \\ a square domain}
  \end{subfigure}
  %
  %
  \begin{subfigure}[b]{0.4 \textwidth}
    {\includegraphics[width=\textwidth]{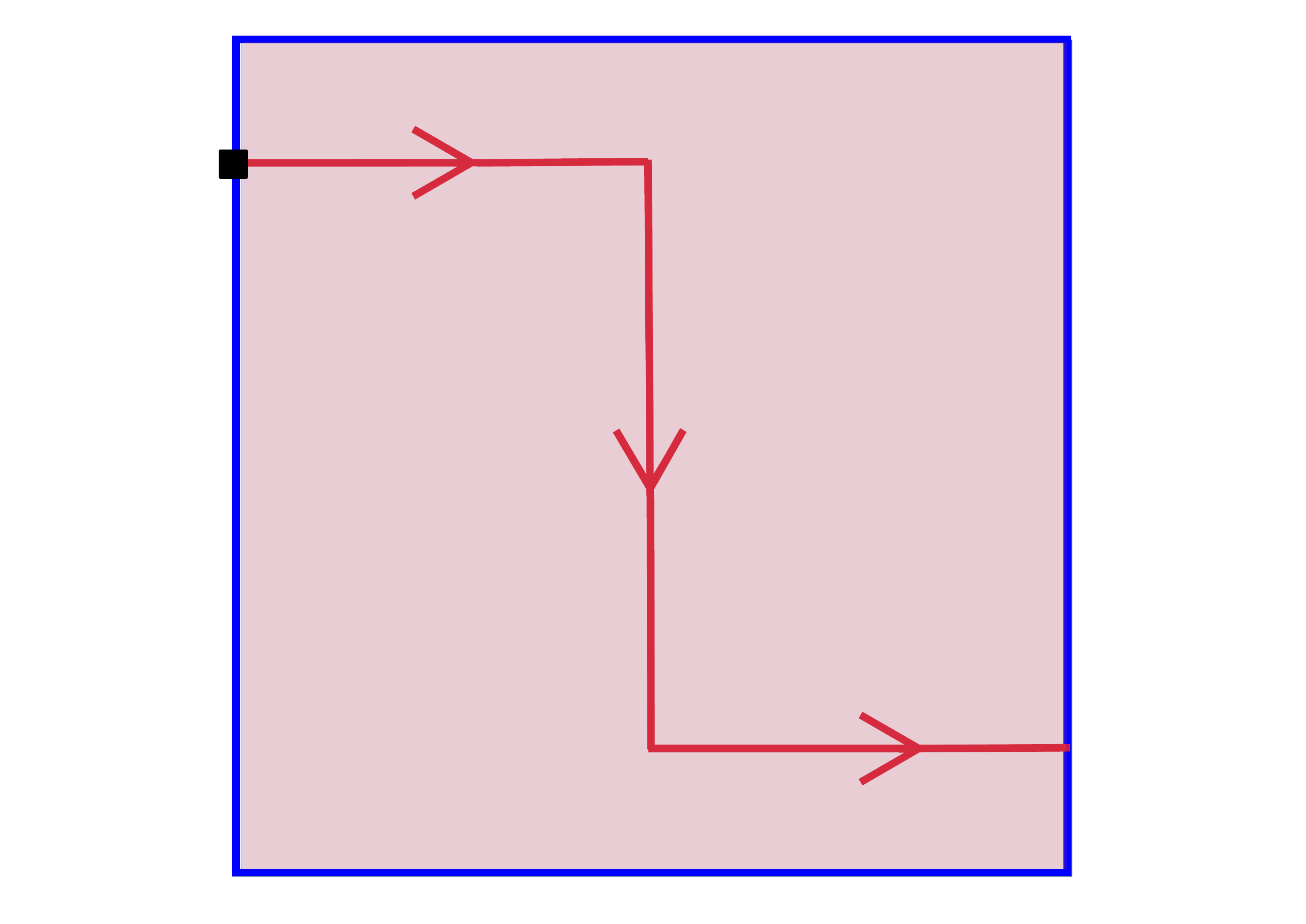}}
    \caption{\tiny Test problem \#4: Stepped channel through a \\ square domain}
  \end{subfigure}
  \begin{subfigure}[b]{0.45 \textwidth}
    {\includegraphics[width=\textwidth]{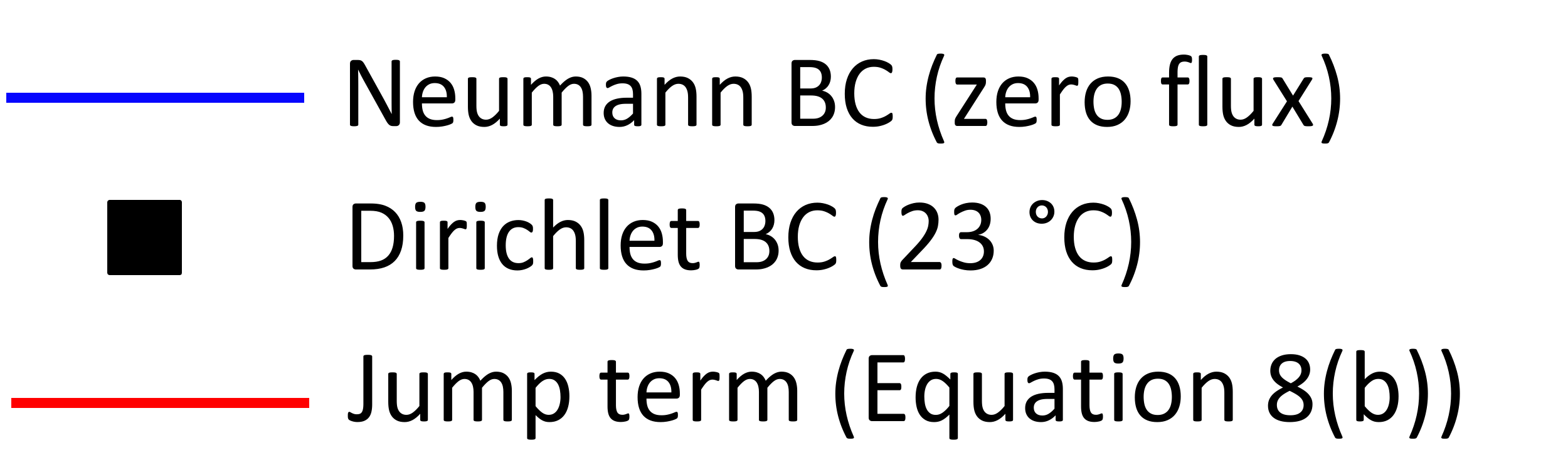}}
    \caption{\tiny Legend for boundary conditions}
  \end{subfigure}
  \caption{The geometries solved using the developed CoolPINNs framework are shown above. The vasculature geometries of straight, curved and zig-zag shape are selected to cover the geometrical configurations of practical interest. All geometries are square with a side length of 100 mm. With origin at the lower left corner of each of the geometries, the vasculature starts at co-ordinate (0, 50) and (0, 25) for geometries A and B, respectively. For geometry C, the vasculature starts at coordinate (0, 50) with sine wave amplitude of 12.5 mm and wavelength of 66.67 mm. For geometry D, the vasculature is defined by co-ordinates (0, 80), (50, 80), (50, 20) and (100, 20). Fluid flow is from left to right through the vasculature for all cases as indicated by the arrows. For the mean surface temperature (MST) study shown in Figs.~\ref{fig:MST_PiNN_No_Rad} and \ref{fig:MST_PiNN_With_Rad}, a slight variation of stepped channel geometry D was used; the coordinates used for MST study were (0, 60), (50, 60), (50, 20) and (100, 20).
  \label{Fig:Geoms_solved}}
\end{figure}

The solutions for all the problems presented in this paper were obtained using the \textsf{DeepXDE} package \citep{lu2021deepxde}. 
The jump term condition, shown by Eq.~\eqref{Eqn:q_jump_condition}, is modeled using the operatorBC command in \textsf{DeepXDE}. 
Note that \textsf{DeepXDE} applies the Dirichlet and Neumann boundary conditions as ``soft constraints'' (i.e., these conditions are not enforced exactly). 
The Dirichlet boundary conditions were applied exactly using a transform function feature available in \textsf{DeepXDE}. The number of uniform collocations points used for the geometries shown in Figs.~\ref{Fig:Geoms_solved}(A), \ref{Fig:Geoms_solved}(B), \ref{Fig:Geoms_solved}(C) and \ref{Fig:Geoms_solved}(D) were 1849, 2025, 1936, and 3942, respectively.

\section{REPRESENTATIVE NUMERICAL RESULTS --- FORWARD PROBLEM}
\label{Sec:S4_PINNS_NR}
In this section, we demonstrate the accuracy of the developed framework using numerical results of test cases. We describe the boundary value problem (BVP) for each test case, summarize the model parameters in the resulting BVPs, and present the CoolPINNs' hyperparameters. We then evaluate the accuracy by comparing the numerical solutions under the proposed methodology with that obtained using the standard Galerkin finite element method.

\subsection{Test case descriptions}

The CoolPINNs modeling framework has been verified by solving four test case problems of a square domain with different shapes of vasculature geometries. 
The vasculature geometries of practical interest include straight, curved and zig-zag with sharp corners. 
The placement of the vasculatures is selected such that the vasculature divides the total domain into symmetric (Fig.~\ref{Fig:Geoms_solved}(A)), non-symmetric (Figs.~\ref{Fig:Geoms_solved}(B) and \ref{Fig:Geoms_solved}(C)), and anti-symmetric (Fig.~\ref{Fig:Geoms_solved}(D)) subdomains. 

A Dirichlet (temperature) boundary condition (i.e., ambient temperature) is specified at a single node at the fluid entry location. 
A Neumann boundary condition is applied on all four boundaries of the square domain by imposing zero flux representing perfectly insulated boundaries. 
The jump boundary condition defined by Eq.~\eqref{Eqn:q_jump_condition} is applied on the entire vasculature. 
The remaining parameters in Eqs.~ \eqref{Eqn:BoE}--\eqref{Eqn:temp_BC_V_in} are defined in Table \ref{tab:params}. 
The fluid enters at the left end and exits from the right end of the vasculature as indicated by the arrows in Fig.~\ref{Fig:Geoms_solved}. 
Note that the coolant material properties defined in Table \ref{tab:params} are for water, and the domain material is a carbon-fiber-reinforced polymer (CFRP). 
An emissivity of 0.95 is used.

Using these test cases, we validated the framework by calculating the error in CoolPINNs solutions with respect to FEM solutions. 
We also verified an important property of constant mean surface temperature of the domain under flow reversal when radiation is excluded. 
To further reinforce the efficacy of our methodology, we solved each of the four problems for wide range of volumetric flow rates of 1 mL/min, 10 mL/min and 100 mL/min.

\begin{table}[h]
    \centering
    \caption{Parameters used in the CoolPINNs numerical simulations. The values are from the experimental literature \citep{devi2021microvascular}.}
\begin{tabular}{lr} \hline 
\textbf{Parameter} & \textbf{Value} \\ \hline
Length $L$ [$\mathrm{mm}$] & 100     \\
Width $B$ [$\mathrm{mm}$]  & 100     \\
Thickness $t$ [$\mathrm{mm}$] & 4.15 \\ 
Emissivity $\epsilon$ [$\mathrm{dimensionless}$] & 0.95 \\ 
Heat transfer coefficient $h_T$ [$\mathrm{W/m^{2}/K}$] & 13.125 \\ 
Applied heater flux $f_0$ [$\mathrm{W/m^2}$] & 1500 \\
Ambient temperature $\vartheta_{\mathrm{amb}}$ [$^{\circ}\mathrm{C}$] & 23 \\ 
Heat capacity of fluid $c_f$ [$\mathrm{J/kg/K}$] & 4800 \\
Density of fluid $\rho_f$ [$\mathrm{kg/m^3}$] & 1000 \\
Volumetric flow rate $V$ [$\mathrm{mL/min}$]  & 1, 10, 100 \\
Thermal conductivity of host CFRP material $K$ [$\mathrm{W/m/K}$] & 2.5247 
\\ \hline 
\end{tabular}
\label{tab:params}
\end{table}

\subsection{Hyperparameter tuning}

As discussed earlier, the CoolPINNs framework attempts to find the solution of the PDE at the collocation points by minimizing the loss function using optimization. 
The enforcement of the losses associated with the PDE (collocation points inside the domain) and boundary conditions (collocation points on the domain boundary) helps reduce the number of possible solutions and improves the framework efficiency. 
However, the non-convexity of the loss function poses some challenges as elaborated below.

In Section \ref{Sec:S1_PINNS_Intro}, we have seen that CoolPINNs framework involves training of deep neural networks (DNNs), objective functions of which are known to be non-convex functions \citep{jain2017non}. 
The non-convexity of the objective or cost function implies presence of multiple local minima. 
This implies that for a given combination of hyperparameters, i.e., number of layers, number of neurons per layer, learning rate etc., the optimization may converge to a specific local minima rather than the actual minima. 
A solution converged to a local minima may not be the most accurate or even can be a wrong solution. Therefore, hyperparameter sweep or tuning is essential to train CoolPINNs for determining the parameters (e.g., learning rate, number of epochs) that enable the framework to converge to the actual minima---needed for getting accurate solutions. 

In the current work, we have performed the hyperparameter sweep for the number of hidden layers, number of neurons per layer, and the learning rate as shown in Table \ref{tab:hyper_sweep}. 
The hyperparameter values used to obtain the results in Figs.~\ref{Fig:Geom1_result}--\ref{Fig:Geom4_result} are shown in the parenthesis in Table \ref{tab:hyper_sweep}.

\begin{table}[h]
    \caption{Hyperparameter sweep performed for the four test cases for forward problems.}
\begin{center}
\begin{tabular}{ |c|c|c|c|c|c| } 
\hline
\textbf{\thead{Test \\ case}} & \textbf{\thead{Coolant \\ flow (mL/min)}} & \textbf{\thead{Hidden \\layers*}} & \textbf{\thead{Neurons \\ per layer*}} & \textbf{\thead{Learning \\ rate*}} & \textbf{\thead{Epochs (Adam \\ + L-BFGS-B**)}}\\
\hline
\multirow{3}{5em}{Straight vasculature at center} & 1 & 3 to 8 (7) & 30 & $10^{-4}$ to $10^{-3}$ ($10^{-3}$) & 10k + 15k\\ 
& 10 & 3 to 7 (4) & 30 & $10^{-4}$ to $10^{-3}$ ($10^{-4}$) & 10k + 15k \\ 
& 100 & 3 to 7 (3) & 30 & $10^{-4}$ to $10^{-3}$ ($8\times10^{-4}$) & 10k + 15k \\ 
\hline
\multirow{3}{5em}{Straight vasculature at quarter} & 1 & 3 to 8 (6) & 30 & $10^{-4}$ to $10^{-3}$ ($7\times10^{-4}$) & 10k + 15k\\ 
& 10 & 3 to 7 (4) & 30 & $10^{-4}$ to $10^{-3}$ ($10^{-4}$) & 10k + 15k \\ 
& 100 & 3 to 7 (3) & 30 & $10^{-4}$ to $10^{-3}$ ($3 \times 10^{-4}$) & 10k + 15k \\ 
\hline
\multirow{3}{5em}{Sine wave vasculature} & 1 & 3 to 8 (3) & 30 & $10^{-4}$ to $10^{-3}$ ($10^{-4}$) & 10k + 15k\\ 
& 10 & 3 to 8 (3) & 30 & $10^{-4}$ to $10^{-3}$ (
$5\times10^{-4}$) & 10k + 15k \\ 
& 100 & 3 to 8 (8) & 30 to 50 (40) & $10^{-4}$ to $10^{-3}$ ($10^{-3}$) & 10k + 15k \\ 
\hline
\multirow{3}{5em}{Stepped vasculature} & 1 & 3 to 10 (10) & 30 to 50 (50) & $10^{-4}$ to $10^{-3}$ ($4\times10^{-4}$) & 25k + 25k\\ 
& 10 & 3 to 10 (5) & 30 to 50 (30) & $10^{-4}$ to $10^{-3}$ ($2\times10^{-4}$) & 10k + 250k \\ 
& 100 & 2 to 10 (9) & 30 to 90 (30) & $10^{-4}$ to $10^{-3}$ ($8\times10^{-4}$) & 25k + 250k \\ 
\hline
\multicolumn{6}{|p{\dimexpr\linewidth-2\tabcolsep-2\arrayrulewidth}|}{*: The numbers in the parenthesis indicate the hyperparameter value used to obtain the results shown in Figs.~\ref{Fig:Geom1_result}--\ref{Fig:Geom4_result}.

**: L-BFGS-B stands for Limited-memory Broyden–Fletcher–Goldfarb–Shanno Boxed algorithm}\\
\hline
\end{tabular}
\end{center}
\label{tab:hyper_sweep}
\end{table}

\subsection{Verification using qualitative properties}
\begin{figure}
    \centering
    \includegraphics[scale=0.5]{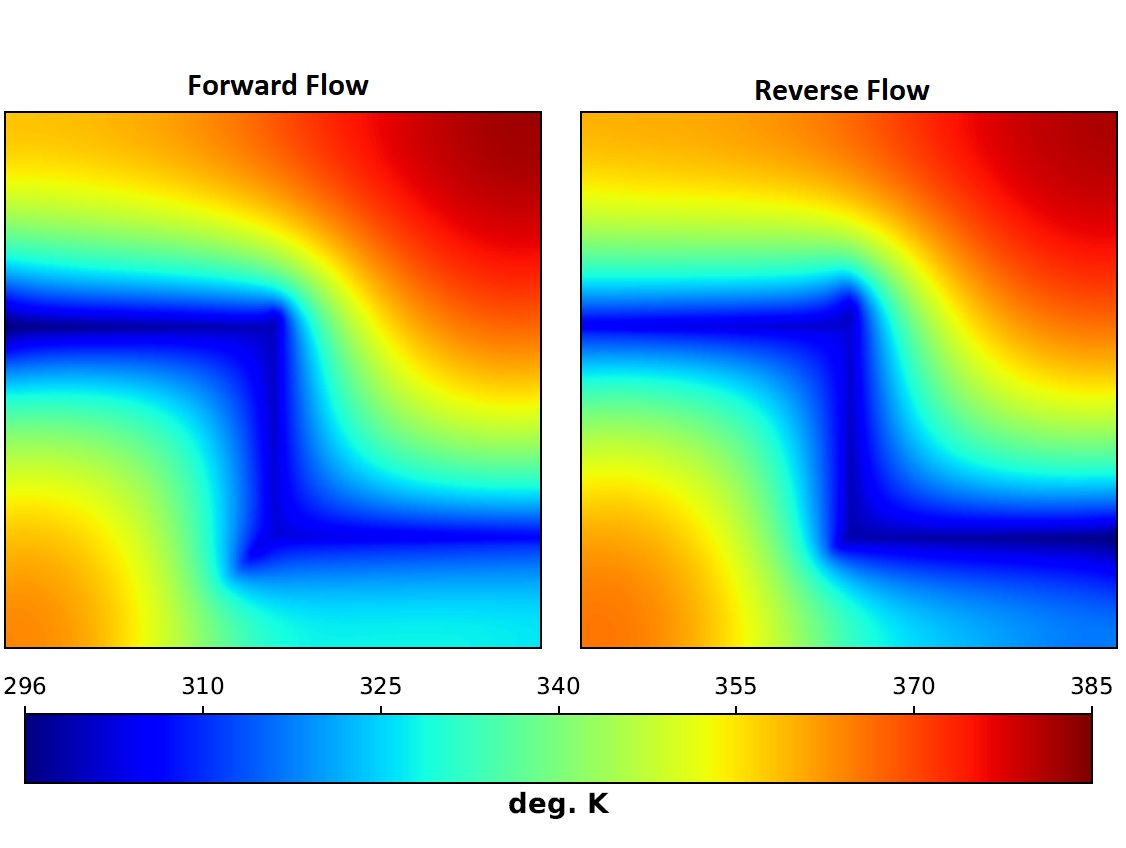}
    \caption{Result obtained for a non-symmetric stepped vasculature with no radiation heat transfer under forward flow (left) and reverse flow (right). The mean surface temperature for the forward and reverse flow cases is 340.60 K and 340.26 K, respectively. 
    The percent difference between these two results is 0.1\%, which shows that the current PINNs architecture can capture the invariance of the mean surface temperature under flow reversal condition.
    \label{fig:MST_PiNN_No_Rad}}
\end{figure}

\begin{figure}
    \centering
    \includegraphics[scale=0.5]{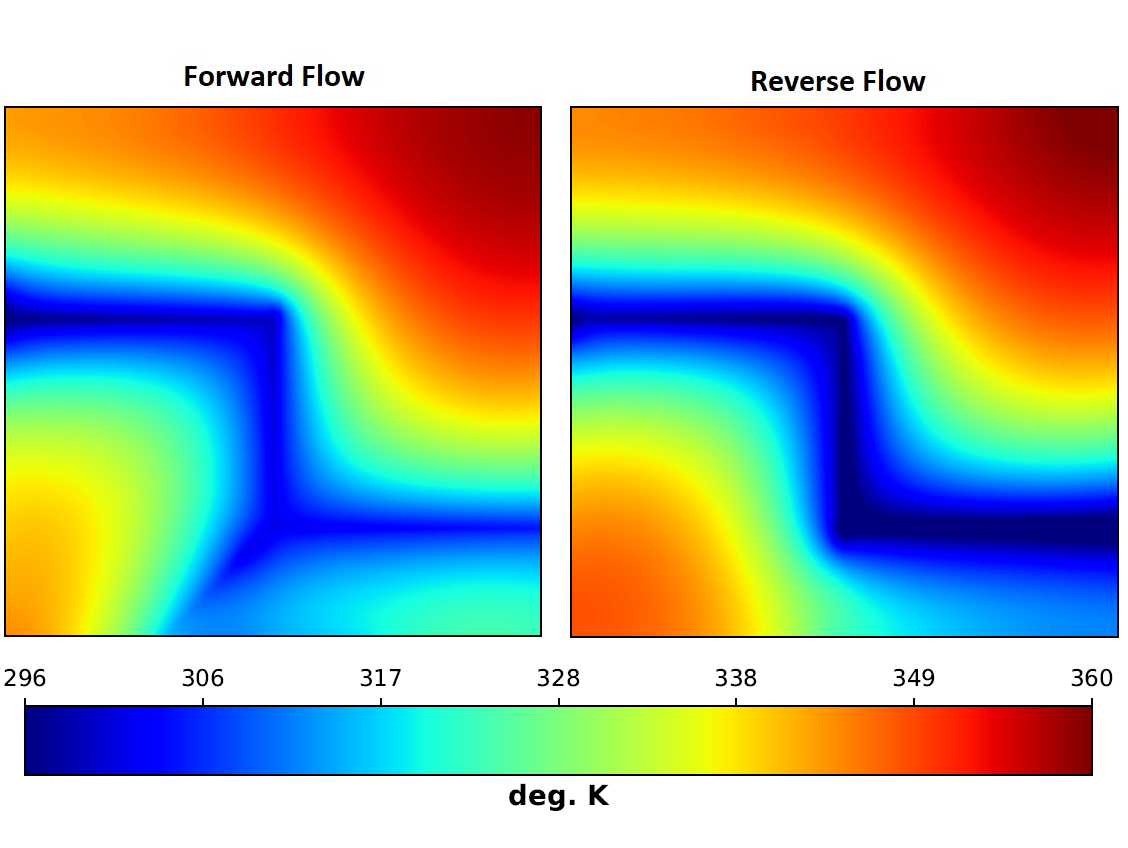}
    \caption{Result obtained for a non-symmetric stepped vasculature with radiation heat transfer included under forward flow (left) and reverse flow (right). The mean surface temperature for the forward and reverse flow cases is 329.62 K and 328.58 K, respectively. 
    When radiation heat transfer is present, the mean surface temperature is not same between forward/reverse flow cases. This is because heat transferred due to radiation between the two cases would be different due to differences in the temperature profiles.
    \label{fig:MST_PiNN_With_Rad}}
\end{figure}

Mean (or average) surface temperature (MST) plays an essential role in the design of slender thermal regulation systems. For instance, MST has been used as one of the metrics to assess performance in studies such as battery packs \citep{pety2017active, WOS:000413131200050},  multifunctional metamaterials \citep{devi2021microvascular},  photovoltaic cells \citep{WOS:000428000700014}, and shape memory alloys \citep{WOS:000320703400007}.

Recently,  \citet{nakshatrala2022Invariance} have shown that MST remains invariant under flow reversal (i.e., swapping the locations of inlet and outlet) if the following conditions are met:
\begin{enumerate}
    \item The domain boundaries are perfectly insulated.
    \item The domain has a constant heat flux source on one face (in other words, the heat generation rate is uniform over the domain).
    \item The thermal load imposed by the vasculature can be idealized as a line load.
    \item Radiation heat transfer is absent.
    \item Conduction and convection are modeled by Fourier model and Newton's law of cooling, respectively.
\end{enumerate}
said differently, the MST under the forward flow condition is the same as the MST under the reverse flow. We will utilize this invariance property to further verify the CoolPINNs modeling framework. Note that, for vasculatures dividing the domain in symmetric sub-domains, the forward and reverse flow problems are essentially the same. 
Therefore, we solved an active-cooling problem in which the embedded vasculature divides the domain into non-symmetric sub-domains. 
For a two-dimensional domain with uniformly distributed collocation points, mean surface temperature is simply an arithmetic average of temperature values at all the collocation points.

Figure \ref{fig:MST_PiNN_No_Rad} shows the results under forward and reverse flow conditions, obtained for a non-symmetric vasculature without radiation. The MST under the forward and reverse flow conditions are 340.60 K and 340.26 K, respectively---the difference between the two results is $0.1 \%$. Therefore, the CoolPINNs architecture respects a fundamental property of vascular-based thermal regulation---MST is invariant under flow reversal.

When radiation heat transfer is also present, the mean surface temperature under flow reversal is not necessarily invariant, but the MST should be close under both flow conditions. Since all other studies in the current work included radiation, we provided the flow reversal results accounting for radiation (see Fig.~\ref{fig:MST_PiNN_With_Rad}). The mean surface temperatures predicted by CoolPINNs under flow reversal conditions are close, demonstrating the framework's effectiveness. 

\subsection{Comparison with FEM solution}
\begin{figure}[h]
  \centering  
    \includegraphics[scale=0.85]{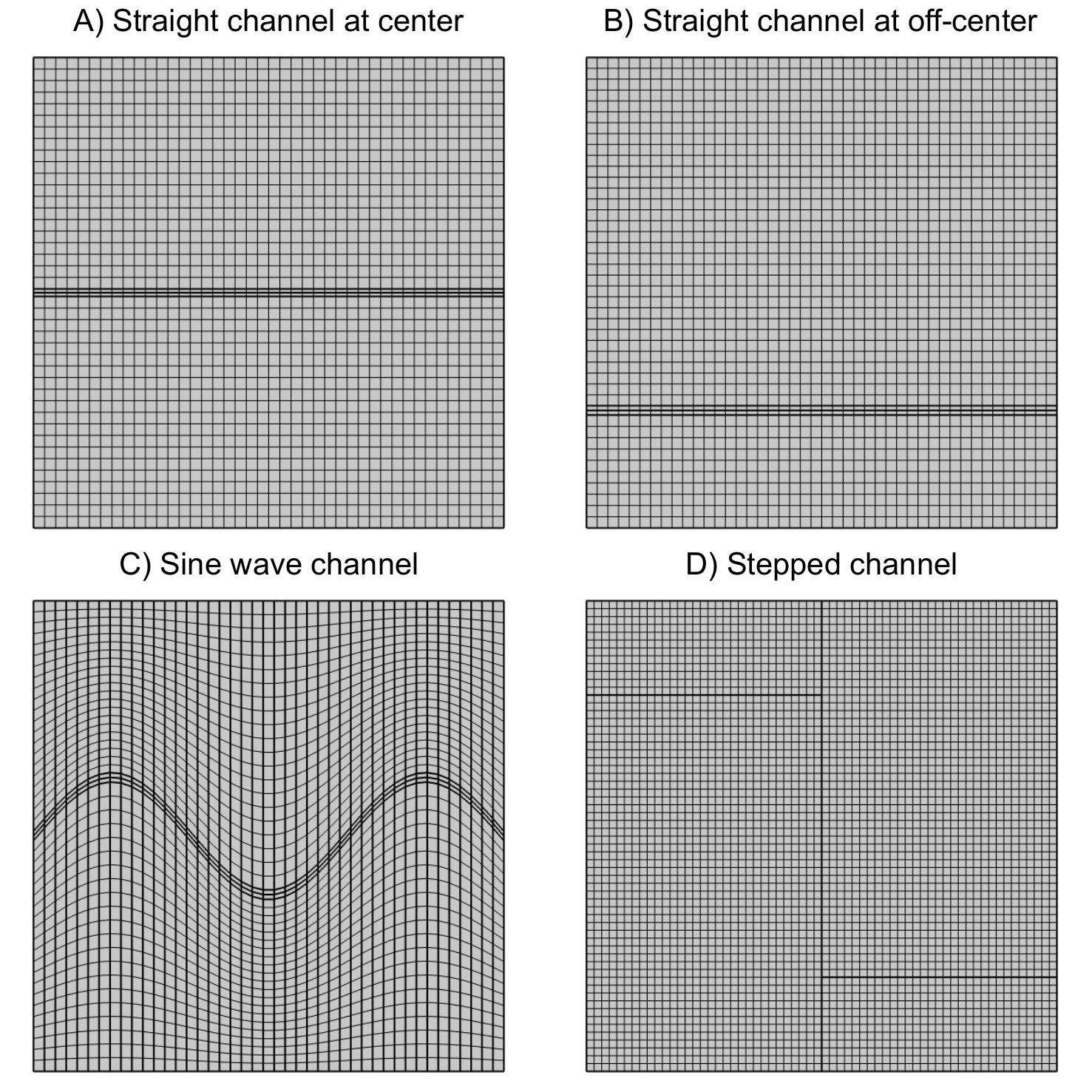}
    \caption{Finite element meshes used in the numerical simulations.}
  \label{Fig:FEM_meshes}
\end{figure}

We have solved the four test cases described earlier using CoolPINNs as well as the finite element method (FEM). For the latter, we have availed the ``weak form'' capability available in \textsf{Comsol}---a multi-physics finite element software package \citep{COMSOL}. Quadrilateral elements with quadratic interpolation were used for spatial discretization under the finite element method. Figure \ref{Fig:FEM_meshes} shows the meshes that were used to solve the four test problems. Appendix \ref{App:Sec_Mesh_refinement} presents the results from a numerical convergence study of finite element solutions
. Below, we compare the two sets of results obtained from PINNs and FEM.

\begin{figure}
  \centering
  \begin{subfigure}[b]{0.9\textwidth}
    {\includegraphics[width=\textwidth]{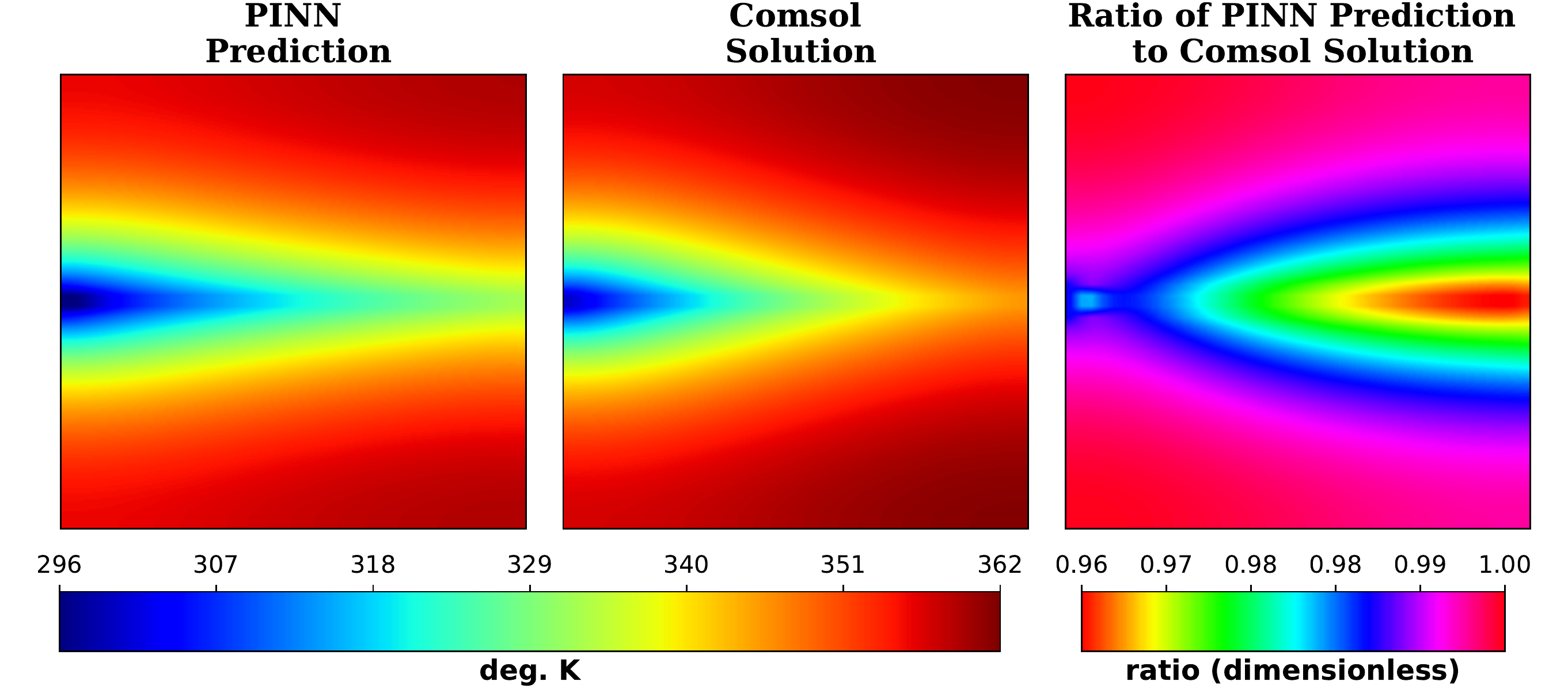}}
  \end{subfigure}
  \par\medskip
  %
  \begin{subfigure}[b]{0.9\textwidth}
    {\includegraphics[width=\textwidth]{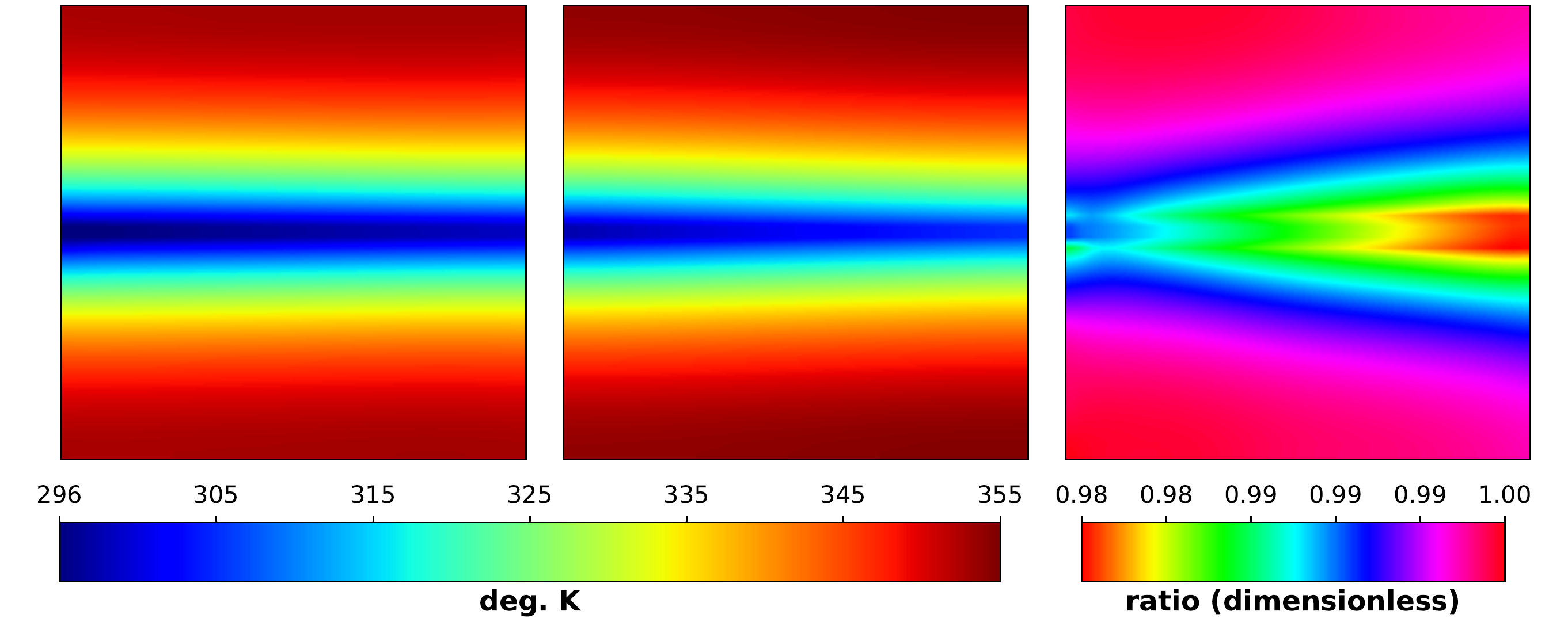}}
  \end{subfigure}
  \par\medskip
  %
  \begin{subfigure}[b]{0.9\textwidth}
    {\includegraphics[width=\textwidth]{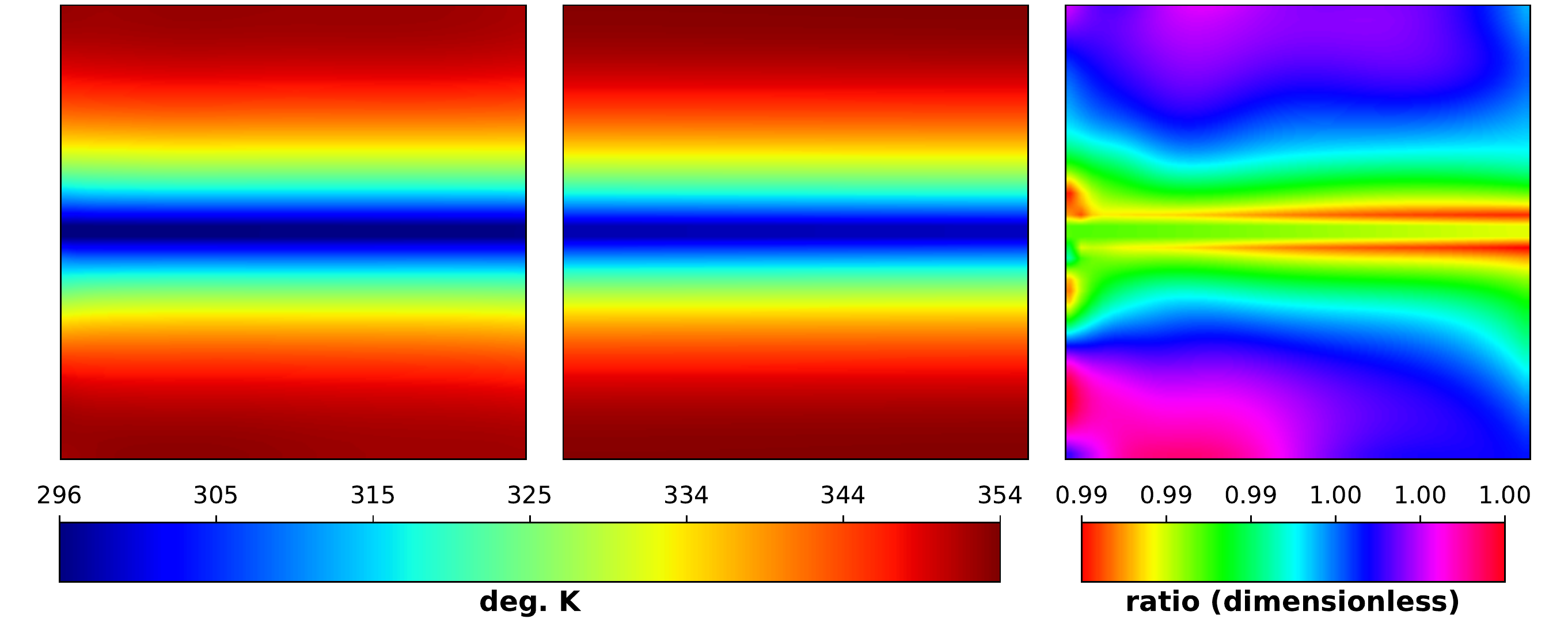}}
  \end{subfigure}
  \caption{Results for test problem \#1 for volumetric flow rates (V) of 1 mL/min, 10 mL/min and 100 mL/min are shown in the top, center and bottom plots. 
  The ratio of CoolPINN prediction to \textsf{Comsol} FEM solution is greater than 0.96, 0.98 and 0.98 for V = 1, 10 and 100 mL/min. The region of relatively high error, indicated by red and yellow colors, for the individual cases is limited to less than 10\% of the domain area from visual observation.}
  \label{Fig:Geom1_result}
\end{figure}

Figure \ref{Fig:Geom1_result} shows the results for a straight vasculature passing through the center of the square domain. 
The PINN prediction comparison with \textsf{Comsol} FEM solutions for this problem for volumetric flow rates of 1 mL/min, 10 mL/min and 100 mL/min are shown in the top, center and bottom plots. 
The PINNs prediction to \textsf{Comsol} solution ratio is greater than 0.96, 0.98 and 0.99 for volumetric flow rates of  1 mL/min, 10 mL/min and 100 mL/min, respectively. 
The PINNs prediction to \textsf{Comsol} solution ratio is lower than 1.0 for all the cases, which indicates that the PINNs is under-predicting the temperatures. 
The region of relatively smaller PINNs-to-\textsf{Comsol} ratio (or larger error) in each plot is indicated by red and yellow colors. The larger error generally appears near the vasculature location. 
The area of these smaller ratio regions for each case is approximately less than 5--10\% of the domain area from visual observation indicating excellent accuracy over majority of the domain area.

\begin{figure}
  \centering
  \begin{subfigure}[b]{0.9\textwidth}
    {\includegraphics[width=\textwidth]{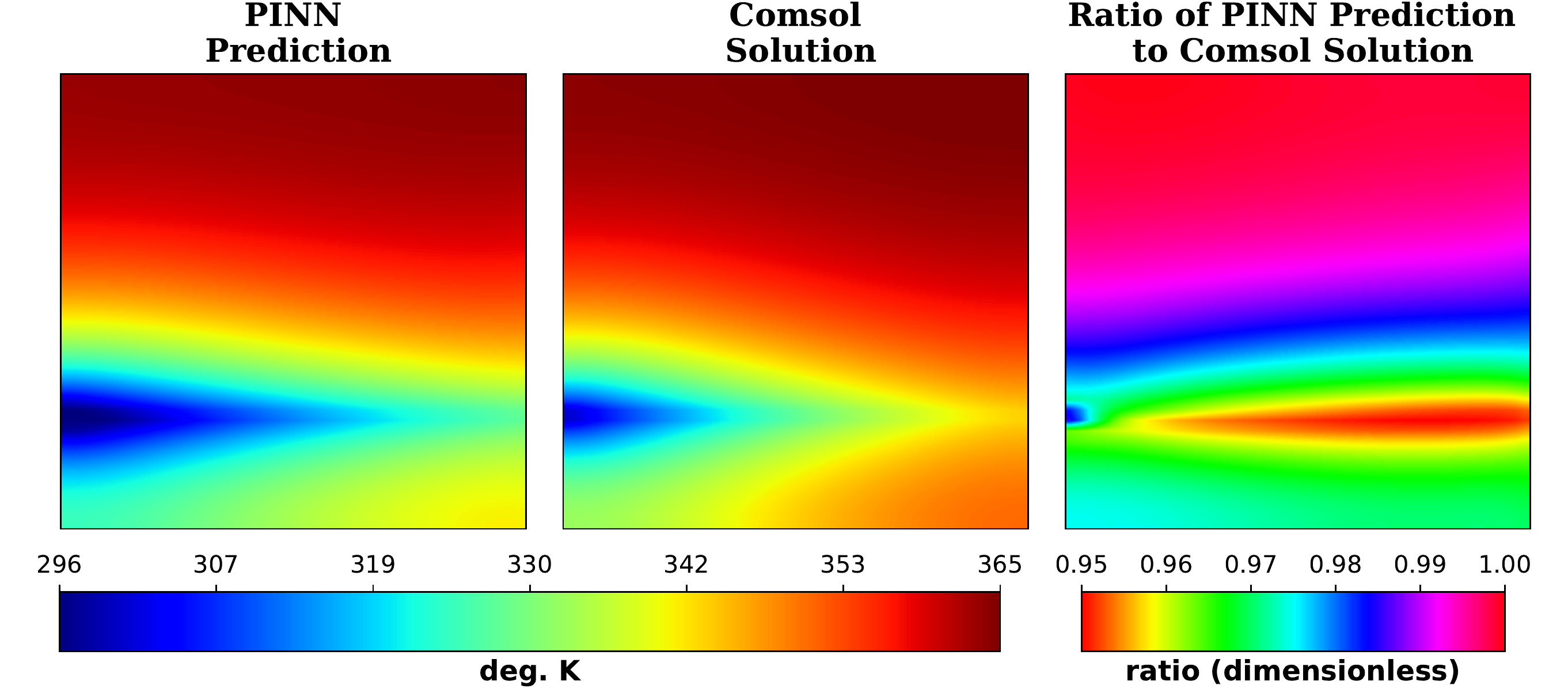}}
  \end{subfigure}
  \par\medskip
  \begin{subfigure}[b]{0.9\textwidth}
    {\includegraphics[width=\textwidth]{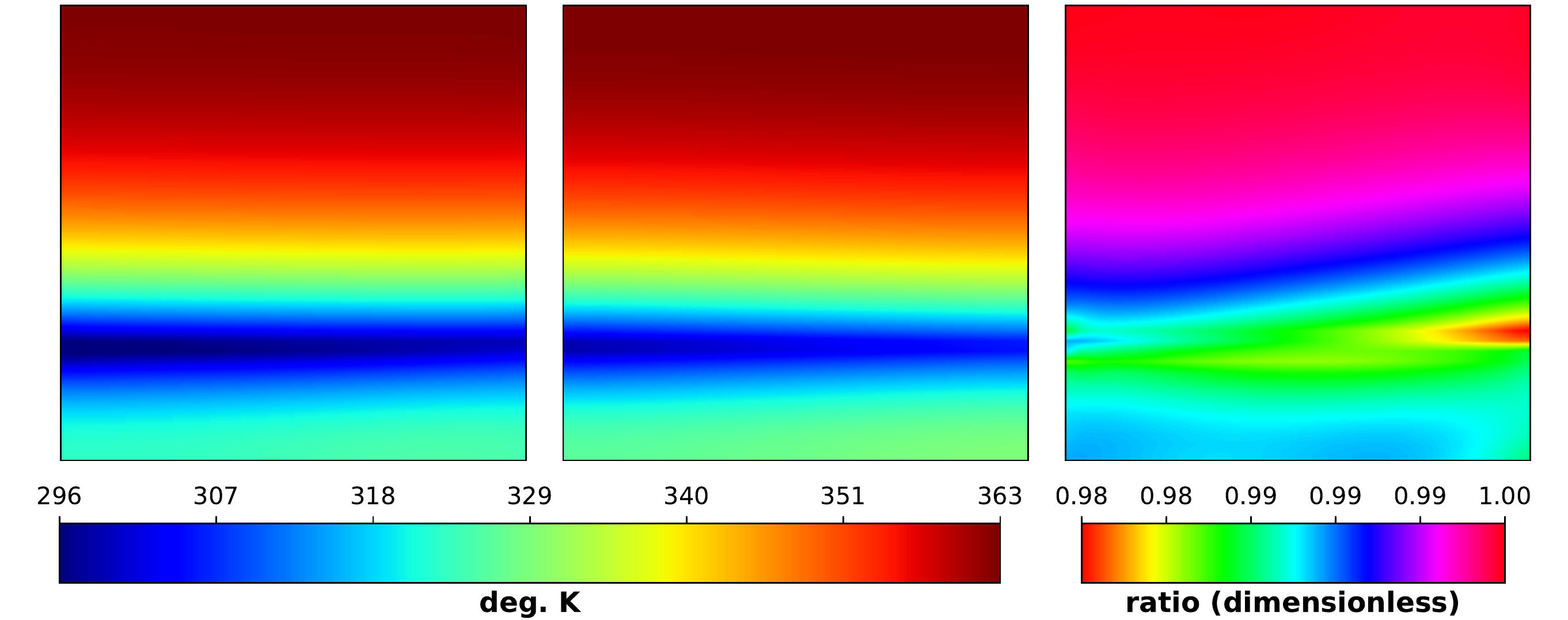}}
  \end{subfigure}
  \par\medskip
  \begin{subfigure}[b]{0.9\textwidth}
    {\includegraphics[width=\textwidth]{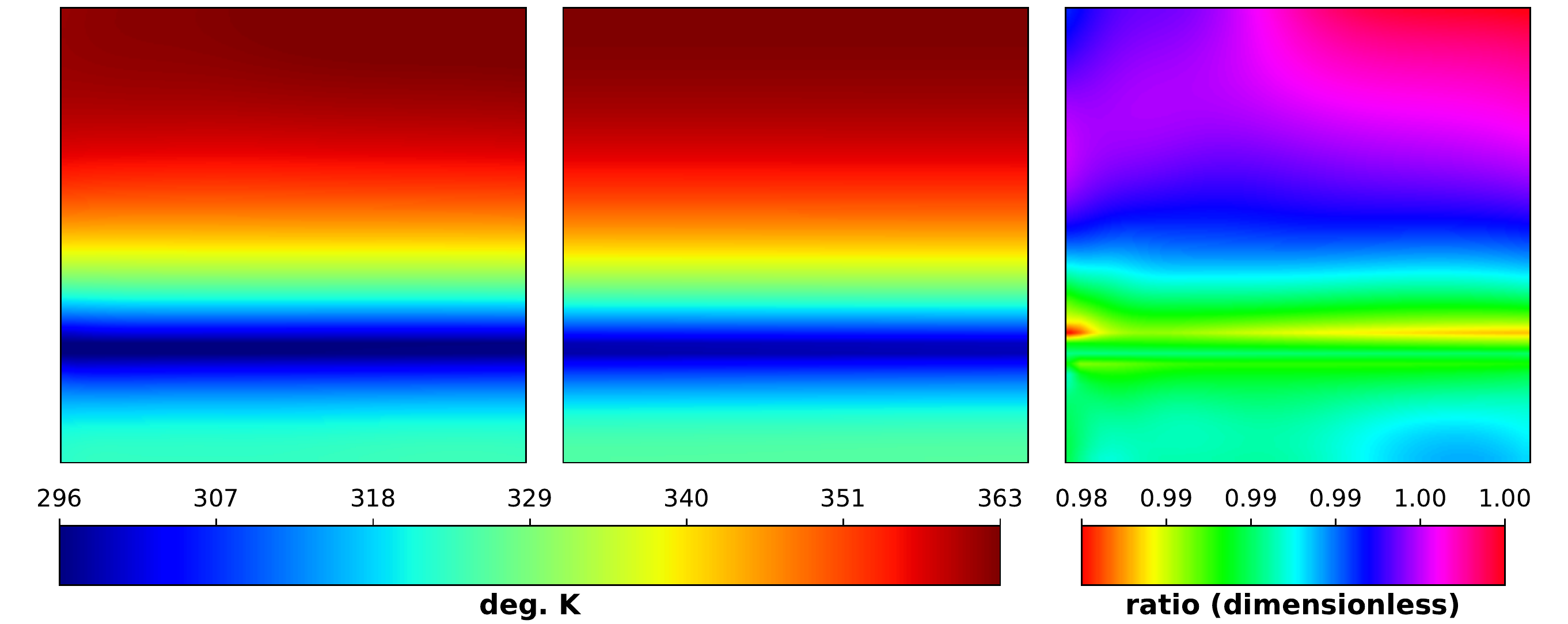}}
  \end{subfigure}
  \caption{Results for test problem \#2 for volumetric flow rates (V) of 1 mL/min, 10 mL/min and 100 mL/min are shown in the top, center and bottom plots. 
  The ratio of PINN prediction to \textsf{Comsol} solution is greater than 0.95, 0.98 and 0.98 for V = 1, 10 and 100 mL/min. 
  The region of relatively high error, indicated by red and yellow colors, for the individual cases is limited to less than 10\% of the domain area from visual observation.}
  \label{Fig:Geom2_result}
\end{figure}

In Fig. \ref{Fig:Geom1_result}, the vasculature divided the domain in a symmetric manner, whereas the next test case in Fig. \ref{Fig:Geom2_result} shows the results for a straight vasculature passing through the quarter height of the square domain dividing the domain in a non-symmetric manner. The PINNs prediction comparison with \textsf{Comsol} FEM solutions for this problem for volumetric flow rates of 1 mL/min, 10 mL/min and 100 mL/min are shown in the top, center and bottom plots. The PINNs prediction to \textsf{Comsol} solution ratio is greater than 0.95 for volumetric flow rates of 1 mL/min and greater than 0.98 for volumetric flow rates of 10 mL/min and 100 mL/min, respectively. The PINNs prediction to \textsf{Comsol} solution ratio is lower than 1.0 for all the cases, which indicates that the PINNs is under-predicting the temperatures. The region of relatively lower PINNs-to-\textsf{Comsol} ratio (or larger error) in each plot are indicated by red and yellow colors. The larger error generally appears near the vasculature location and the area of these lower ratio regions for each case is approximately less than 5-10\% of the domain area from visual observation, indicating excellent accuracy over the majority of the domain area.

\begin{figure}
  \centering
  \begin{subfigure}[b]{0.86\textwidth}
    {\includegraphics[width=\textwidth]{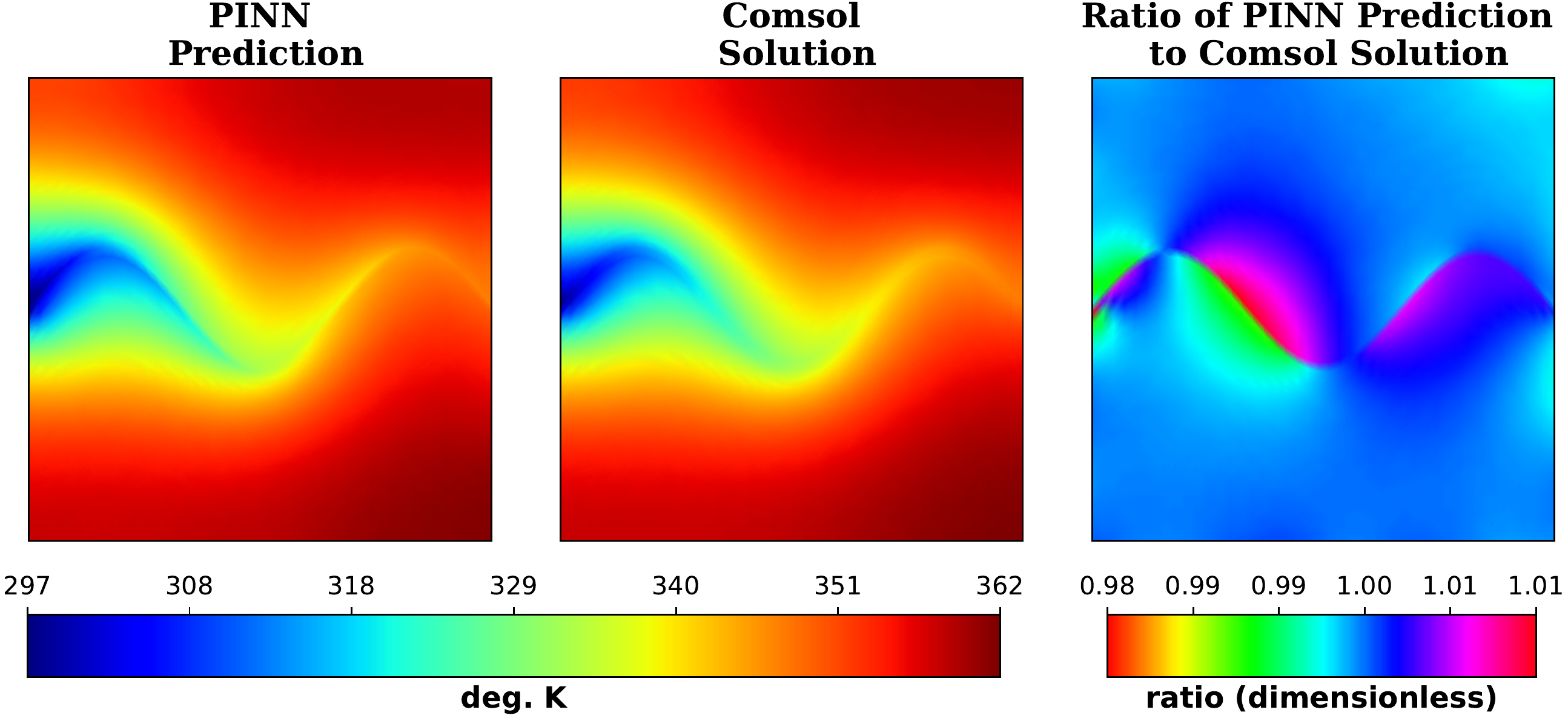}}
  \end{subfigure}
  \par\medskip
  \begin{subfigure}[b]{0.86\textwidth}
    {\includegraphics[width=\textwidth]{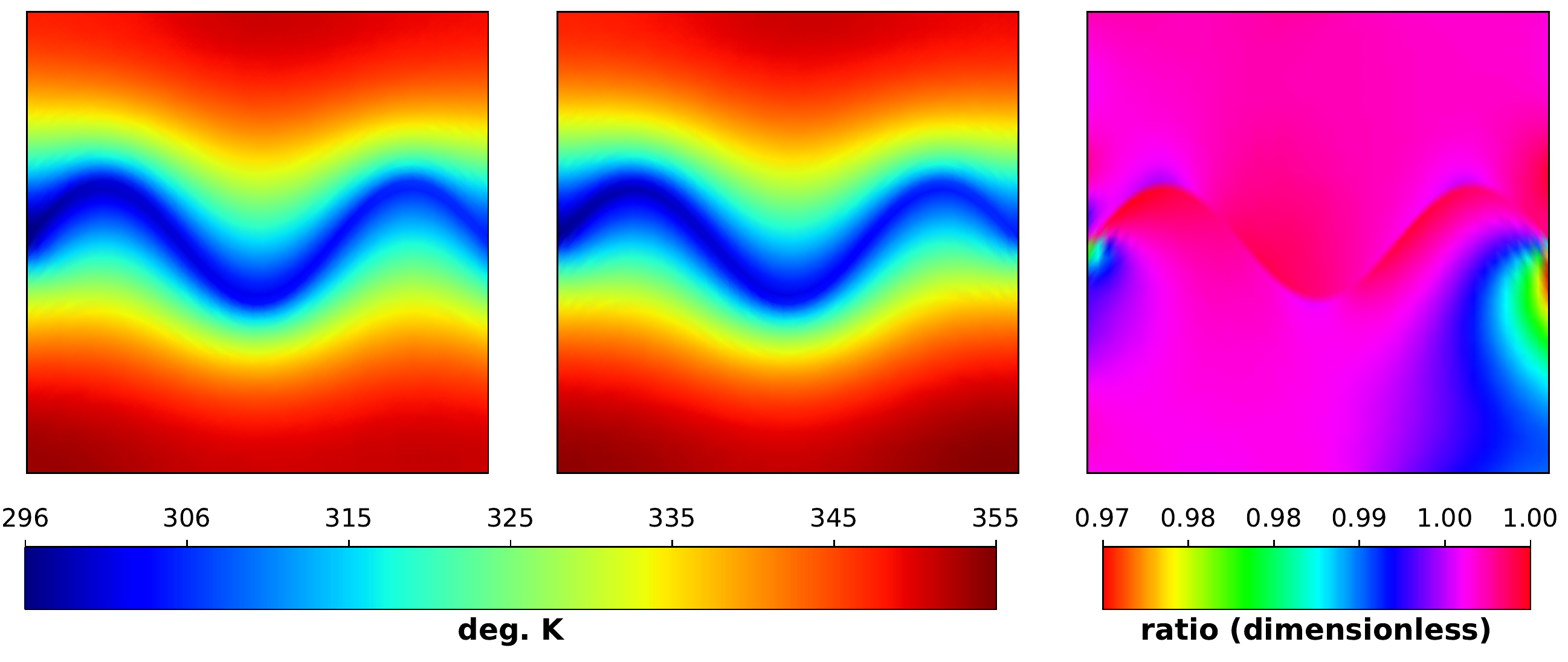}}
  \end{subfigure}
  \par\medskip
  \begin{subfigure}[b]{0.86\textwidth}
    {\includegraphics[width=\textwidth]{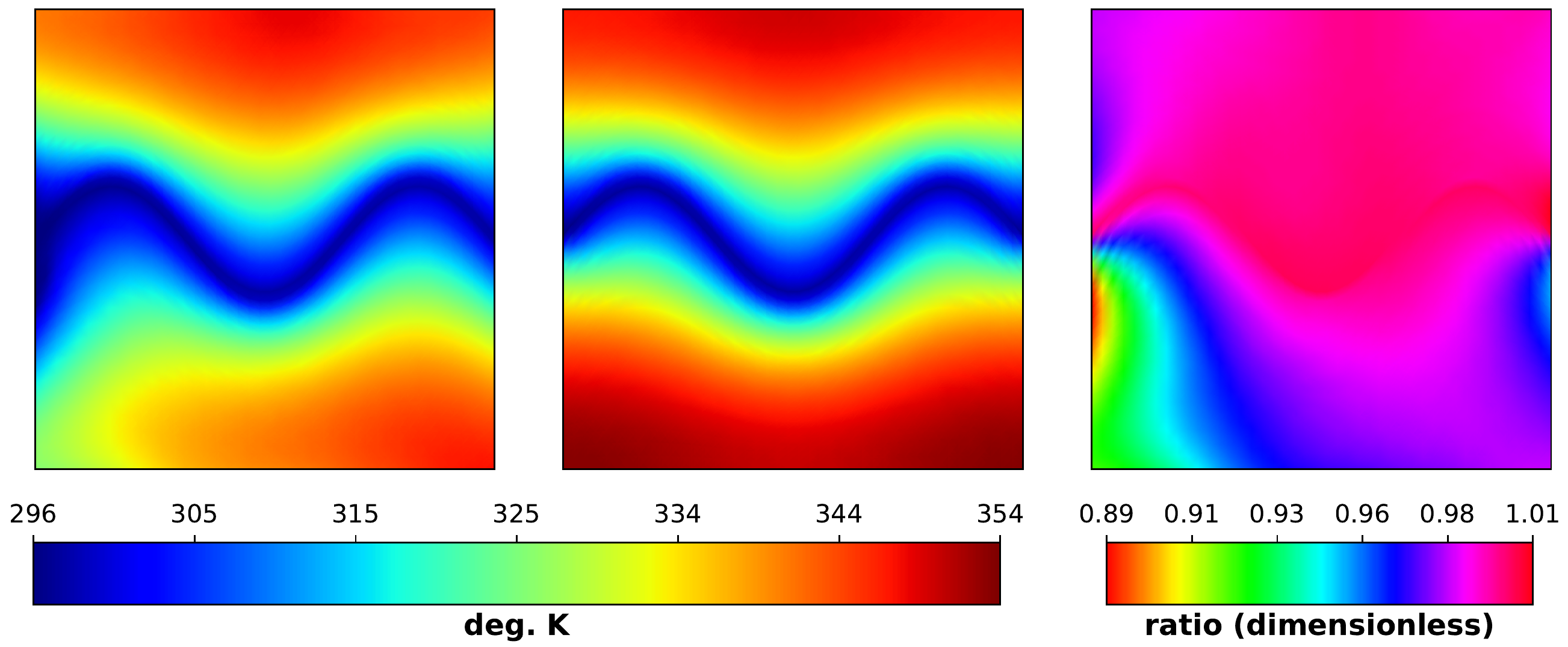}}
  \end{subfigure}
  \caption{Results for test problem \#3 for volumetric flow rates (V) of 1 mL/min, 10 mL/min and 100 mL/min are shown in the top, center and bottom plots. 
  The ratio of PINN prediction to \textsf{Comsol} solution is greater than 0.98, 0.97 and 0.89 for V = 1, 10 and 100 mL/min. 
  The region of relatively high error, indicated by red and yellow colors, for the individual cases is limited to less than 10\% of the domain area from visual observation.}
  \label{Fig:Geom3_result}
\end{figure}

\begin{figure}
  \centering
  \begin{subfigure}[b]{0.9\textwidth}
    {\includegraphics[width=\textwidth]{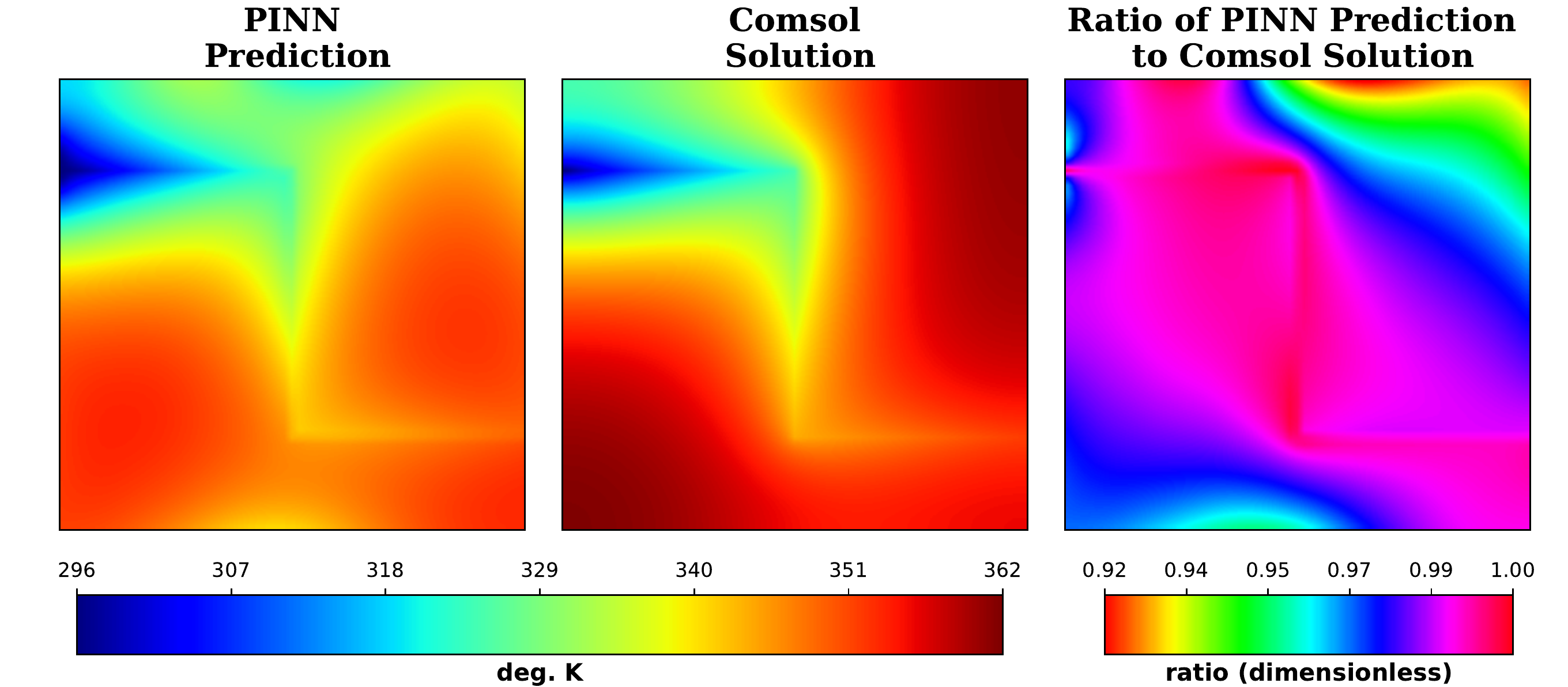}}
  \end{subfigure}
  \par\medskip
  \begin{subfigure}[b]{0.9\textwidth}
    {\includegraphics[width=\textwidth]{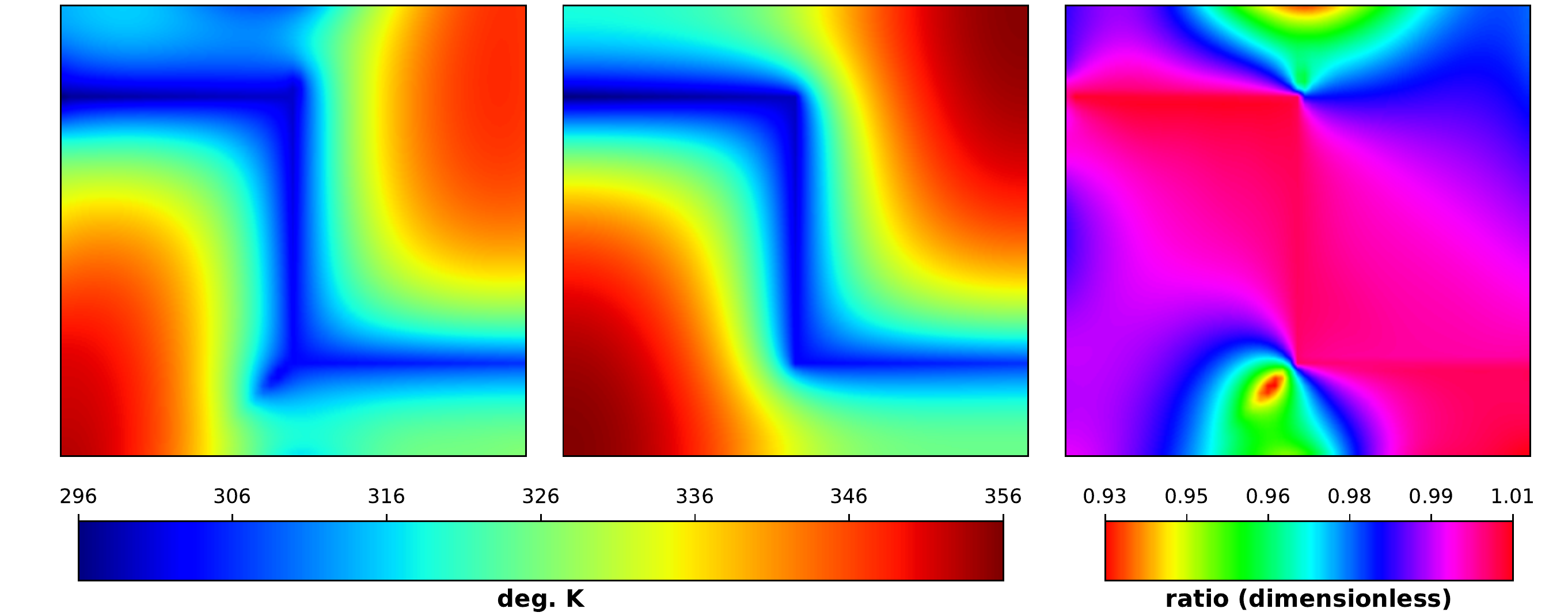}}
  \end{subfigure}
  \par\medskip
  \begin{subfigure}[b]{0.9\textwidth}
    {\includegraphics[width=\textwidth]{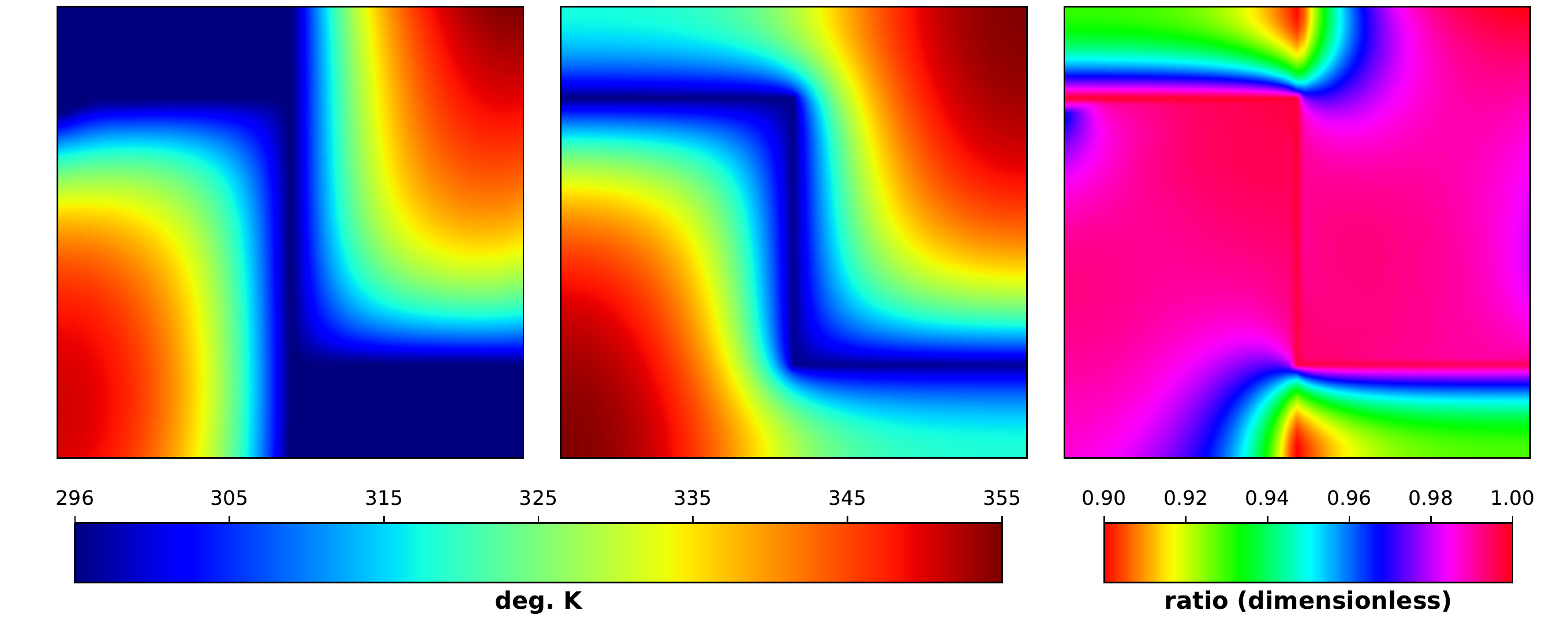}}
  \end{subfigure}
  \caption{Results for test problem \#4 for volumetric flow rates (V) of 1 mL/min, 10 mL/min and 100 mL/min are shown in the top, center and bottom plots. The ratio of PINN prediction to \textsf{Comsol} solution is greater than 0.92, 0.93 and 0.90 for V = 1, 10 and 100 mL/min. The region of relatively high error, indicated by red and yellow colors, is limited to less than 10\% of the domain area from visual observation.}
  \label{Fig:Geom4_result}
\end{figure}

The vasculatures analyzed in Figs. \ref{Fig:Geom1_result} and \ref{Fig:Geom2_result} are straight; however, in real-life for efficient heat transfer the vasculature shapes are often complex. Figures \ref{Fig:Geom3_result} and \ref{Fig:Geom4_result} demonstrate the efficacy of the framework to solve moderately complex vasculature shapes of sine wave and a single step. Fig. \ref{Fig:Geom3_result} shows the results for a sine-wave vasculature originating at the center of the left side of the square domain. The PINNs prediction comparison with \textsf{Comsol} FEM solutions for this problem for volumetric flow rates of 1 mL/min, 10 mL/min and 100 mL/min are shown in the top, center and bottom plots. Compared with the \textsf{Comsol} FEM solutions, the PINNs prediction to \textsf{Comsol} solution ratio is greater than 0.98, 0.97 and 0.89 for the volumetric flow rates of 1, 10 and 100 mL/min, respectively. The ratio the PINNs prediction to \textsf{Comsol} solution ranges between 0.87 to 1.01 for the three cases, which indicates that the PINNs is generally under-predicting the temperatures. The region of relatively lower PINNs-to-\textsf{Comsol} ratio (or larger error) in each plot are indicated by red and yellow colors. For the volumetric flow rates of 1 and 10 mL/min, the larger error appears near the coolant entry and exit locations on the vasculature and the majority of the remaining domain is covered by blue, purple and magenta colors indicating extremely high accuracy of PINNs prediction. For the coolant flow rate of 100 mL/min, the Neumann boundary condition near the lower left side is not accurately enforced. Note that, DeepXDE package applies Dirichlet and Neumann boundary conditions as a soft constraint, and the Neumann boundary conditions cannot be enforced as a hard constraint without custom programming.

Figure \ref{Fig:Geom4_result} shows the results for a single step vasculature. The PINNs prediction comparison with \textsf{Comsol} FEM solutions for this problem for volumetric flow rates of 1 mL/min, 10 mL/min and 100 mL/min are shown in the top, center and bottom plots. The ratio the PINNs prediction to \textsf{Comsol} solution is greater than 0.92, 0.93 and 0.90 for volumetric flow rates of 1, 10 and 100 mL/min, respectively. The ratio the PINNs prediction to \textsf{Comsol} solution ranges between 0.90 to 1.01 for the three cases, and therefore the PINNs is generally under-predicting the temperatures. The region of relatively lower PINNs-to-\textsf{Comsol} ratio (or larger error) in each plot are indicated by red and yellow colors. For all the three volumetric flow rates, the larger error appears near the edges of the domain indicating that the Neumann boundary condition in these regions is not accurately enforced accurately. In the DeepXDE package, the Dirichlet and Neumann boundary conditions are enforced as soft constraints, and the Neumann boundary condition cannot be enforced as a hard constraint unlike the Dirichlet boundary condition. The majority of the remainder domain is covered by blue, purple and magenta colors indicating extremely high accuracy of PINNs prediction.

\section{REPRESENTATIVE NUMERICAL RESULTS --- INVERSE PROBLEM}
\label{Sec:S5_PINNS_INV}

The same test cases used for forward problems, shown in Fig.~\ref{Fig:Geoms_solved}, were used to solve the inverse problems. 
The parameters used for the inverse problems were the same as forward problem (from Table \ref{tab:params}) with the exception of the thermal conductivity, which was the parameter predicted using the inverse solution.
We took the coolant's  volumetric flow rate of 10 mL/min. 
The number of random collocations points used for the geometries shown in Figs.~\ref{Fig:Geoms_solved}(A), \ref{Fig:Geoms_solved}(B), \ref{Fig:Geoms_solved}(C) and \ref{Fig:Geoms_solved}(D) were 1804, 1754, 2541, and 2102, respectively. 
Similar to forward problems, a hyperparameter sweep for training the model was performed to determine the model providing the most accurate prediction.
The hyperparameter sweep study is summarized in Table \ref{tab:hyper_sweep_inverse}.
To demonstrate the effectiveness of the framework to solve the inverse problems, two different initial values for thermal conductivity, 1.0 W/m/K and 6.25 W/m/K, were used compared to the ground truth of 2.5247 W/m/K. 
The temperature solutions obtained from the forward problem were used for solving the inverse problems. 
Since the temperature solution from the forward problem has some amount of error compared to the actual solution, temperatures used as input to the inverse problems can be considered as noisy data; e.g., temperatures obtained from a real-life scenario like thermal imagery.

\begin{figure}
    \centering
    \includegraphics[scale=0.435]{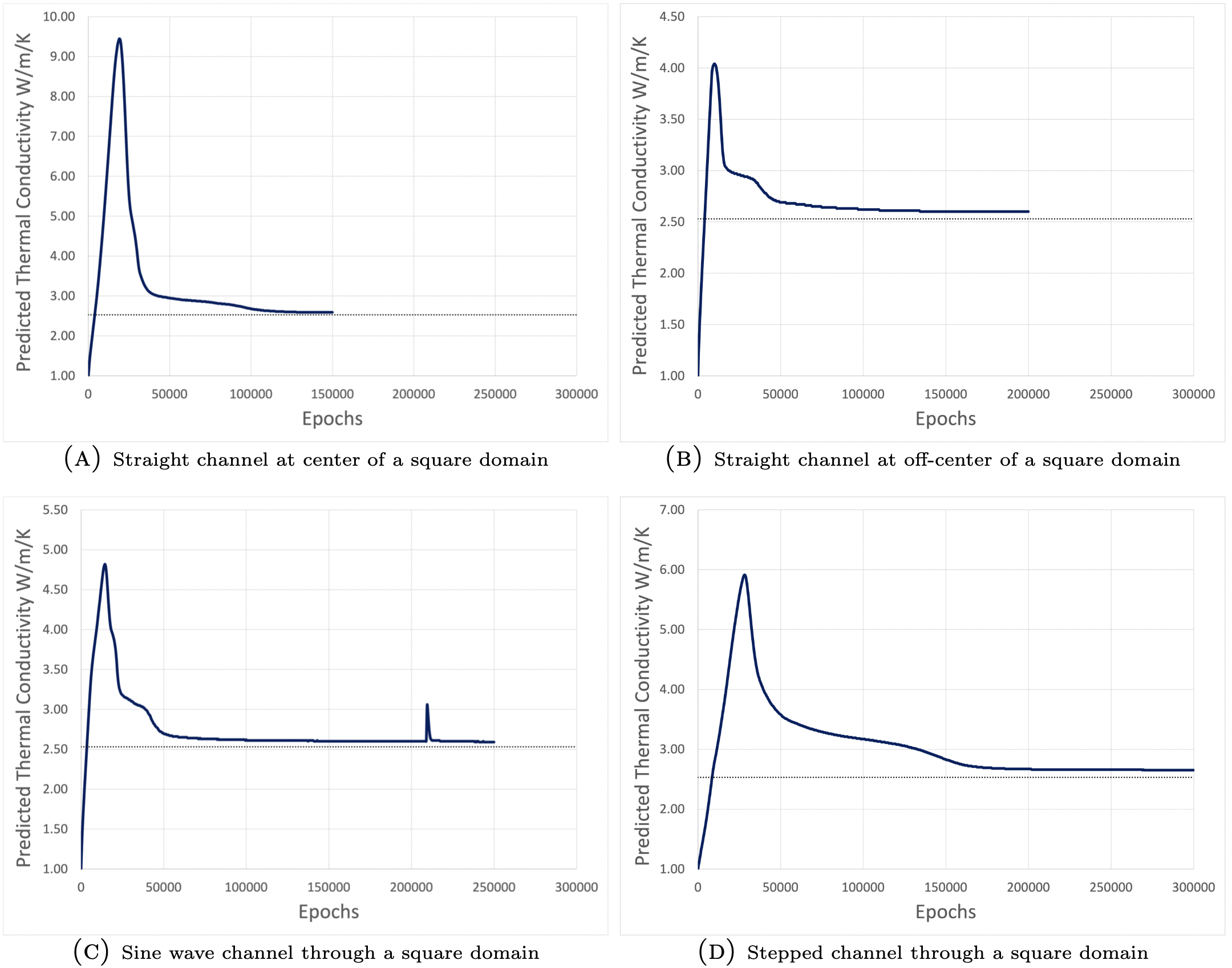}
    \caption{Variation of predicted thermal conductivity with respect to the number of epochs used to train the model for the initial thermal conductivity of 1 W/m/K. The dotted black line in all plots represents the ground truth of 2.5247 W/m/K. The final predictions for the four test cases are 2.59, 2.64, 2.60, and 2.67 W/m/K.
    \label{Fig:Inverse_prob_solved_Kprst=1}}
\end{figure}

\begin{figure}
    \centering
    \includegraphics[scale=0.435]{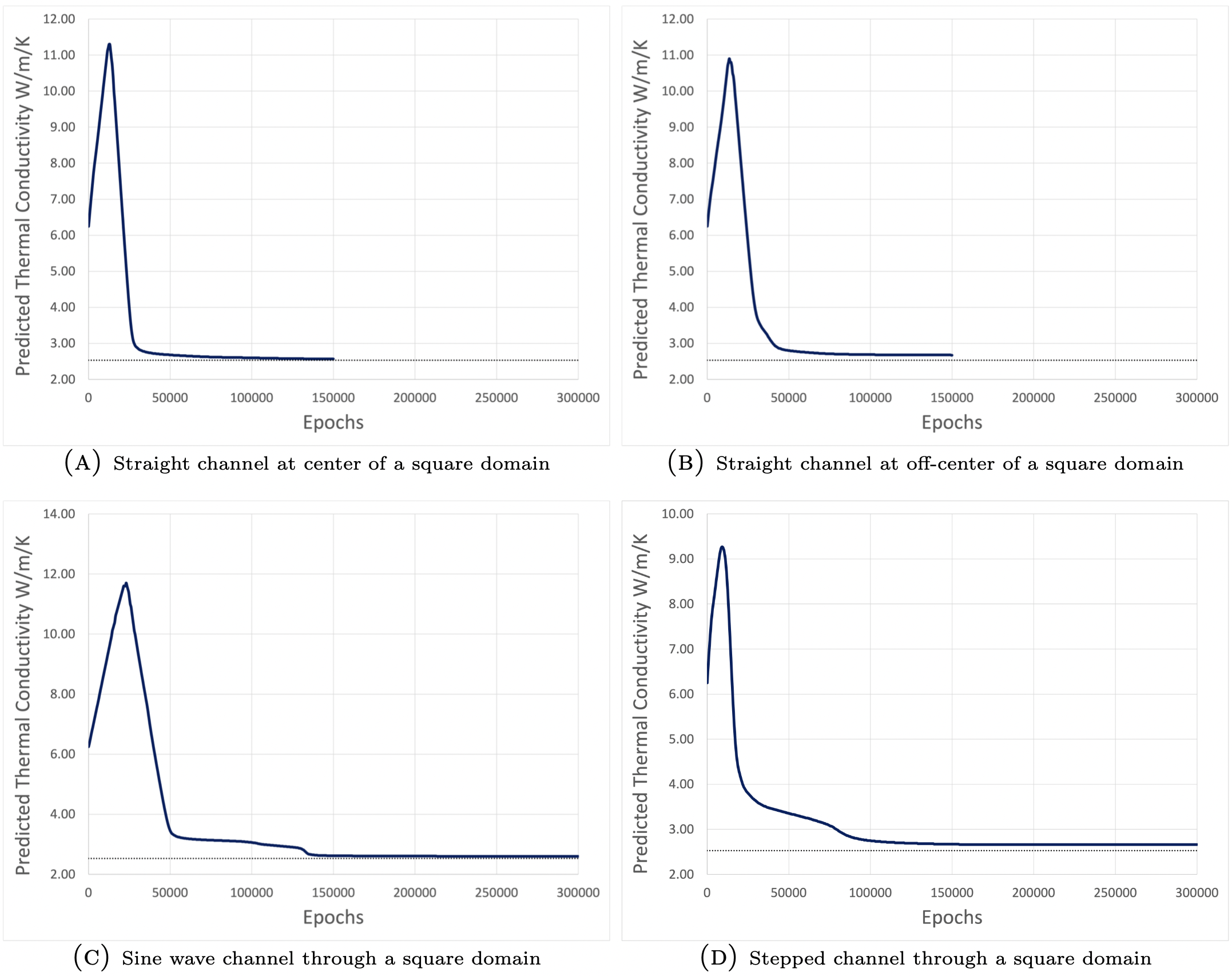}
    \caption{Variation of predicted thermal conductivity with respect to the number of epochs used to train the model for the initial thermal conductivity of 6.25 W/m/K. The dotted black line in all plots represents the ground truth of 2.5247 W/m/K. The final predictions for the four test cases are 2.57, 2.67, 2.60, and 2.66 W/m/K.
    \label{Fig:Inverse_prob_solved_Kprst=6.25}}
\end{figure}

\begin{table}[h]
    \caption{Hyperparameter sweep performed for the four test cases for inverse problems.}
\begin{center}
\begin{tabular}{ |c|c|c|c|c|c| } 
\hline
\textbf{\thead{Test \\ case}} & \textbf{\thead{Initial K \\ (W/m/K)}} & \textbf{\thead{Hidden \\layers*}} & \textbf{\thead{Neurons \\ per layer}} & \textbf{\thead{Learning \\ rate*}} & \textbf{\thead{Epochs (Adam)}}\\
\hline
\multirow{3}{5em}{Straight vasculature at center} 
& 1 & 2 to 3 (2) & 30 & $10^{-4}$ to $10^{-3}$ ($9 \times 10^{-4}$) & 150k\\ 
& 6.25 & 2 to 3 (3) & 30 & $10^{-4}$ to $10^{-3}$ ($6 \times 10^{-4}$) & 150k\\ 
&      &            &    &                     &\\
\hline
\multirow{3}{5em}{Straight vasculature at quarter} 
& 1 & 3 to 4 (2) & 30 to 40 & $10^{-4}$ to $10^{-3}$ ($7 \times 10^{-4}$) & 200k\\ 
& 6.25 & 3 to 4 (3) & 30 & $10^{-4}$ to $10^{-3}$ ($5 \times 10^{-4}$) & 150k\\ 
&      &            &    &                     &\\
\hline
\multirow{3}{5em}{Sine wave vasculature} 
& 1 & 4 to 5 (4) & 40 & $10^{-4}$ to $10^{-3}$ ($10^{-3}$) & 250k\\ 
& 6.25 & 4 to 5 & 40 & $10^{-4}$ to $10^{-3}$ ($3 \times 10^{-4}$) & 300k\\ 
&      &            &    &                     &\\
\hline
\multirow{3}{5em}{Stepped vasculature} 
& 1 & 3 to 5 (3) & 30 & $10^{-4}$ to $10^{-3}$ ($3 \times 10^{-4}$) & 300k\\ 
& 6.25 & 3 to 5 (3) & 30 & $10^{-4}$ to $10^{-3}$ ($7 \times 10^{-4}$) & 300k\\ 
&      &            &    &                     &\\
\hline
\multicolumn{6}{|p{\dimexpr\linewidth-2\tabcolsep-2\arrayrulewidth}|}{*The numbers in the parenthesis are the hyperparameter values used to obtain the results shown in Figs.~\ref{Fig:Inverse_prob_solved_Kprst=1} and \ref{Fig:Inverse_prob_solved_Kprst=6.25}.}\\
\hline
\end{tabular}
\end{center}
\label{tab:hyper_sweep_inverse}
\end{table}

Figures \ref{Fig:Inverse_prob_solved_Kprst=1} and \ref{Fig:Inverse_prob_solved_Kprst=6.25} show the thermal conductivity predictions using CoolPINNs framework for initial values for thermal conductivity values of 1.0 W/m/K and 6.25 W/m/K. 
The dotted black line in all plots represents the ground truth of 2.5247 W/m/K. 
The predictions converge to the ground truth for both initial values of thermal conductivity that are higher and lower than the ground truth.
Therefore, the predictions are insensitive to the initial value of the variable to be predicted, which is a major advantage of the CoolPINNs framework.
The percent error with respect to the ground truth is between 1.8\% to 5.8\% for all the cases solved, showing very good accuracy.
As the plots indicate, the accuracy can further be improved by training the model for larger number of epochs.

\section{CLOSURE}
\label{Sec:S6_PINNS_Closure}
We presented CoolPINNs framework based on physics-informed neural networks (PINNs) to model thermal regulation in thin domains with vasculatures. 
The framework solves forward and inverse problems for a reduced-order, non-linear heat transfer PDE with a jump term in thermal flux in a robust manner. 
The framework is a meshless method and provides fast forecast.

The CoolPINNs framework for forward problems was verified by solving four problems with different vasculature geometries, each with three widely varying coolant velocities. 
The solutions obtained using the PINNs methodology were compared with FEM solutions to demonstrate excellent accuracy.
The robustness of the framework for inverse problems was demonstrated by predicting thermal conductivity of the panel material using noisy temperature data for all four problems.
The framework provided accurate predictions using two different initial values of thermal conductivity; higher and lower than the ground truth, indicating that the framework is insensitive to initial value of the variable to be predicted. 
A broad sweep of hyperparameters was performed for both forward and inverse problems. 

Although the PINNs framework is developed for a reduced-order mathematical model---pertinent to thermal regulation in several emerging technologies (e.g., high-powered antennas), the approach equally applies to three-dimensional problems and to other thermal regulation models. Additionally, the framework for inverse problems can be easily extended to predict magnitude of the heat source, boundary flux and bulk radiation properties by making only a few changes in the code.
Solution of these problems are very critical for applications like predicting material degradation of heat shield panels and monitoring radiation properties of the reflecting surfaces of cryogenic panels.

We envision that the CoolPINNs framework can greatly benefit the scientific community to solve inverse problems that are difficult and extremely time-consuming to handle using the traditional methods. Compared to grid-based methods (e.g., finite element method and finite volume method), PINN-based methods do not need computational grids and thus can handle complex geometries more easily.
Solution of the specific inverse problem, e.g., parameter identification, can be used to automate the process of informed decision making. For extreme environment applications like Lunar and Martian missions, vasculature panel's surface temperature profile can be periodically captured to inversely calculate the thermal conductivity. The changes in the thermal conductivity can be monitored to identify the onset of material degradation. The rate of degradation can be decelerated by taking appropriate actions like increasing coolant circulation. The entire sequence of capturing temperature profile to increase the coolant circulation can be used to develop an autonomous system for informed decision making.

\appendix

\section{Mesh refinement studies for finite element solutions}
\label{App:Sec_Mesh_refinement}

\begin{figure}
    \centering
    \includegraphics{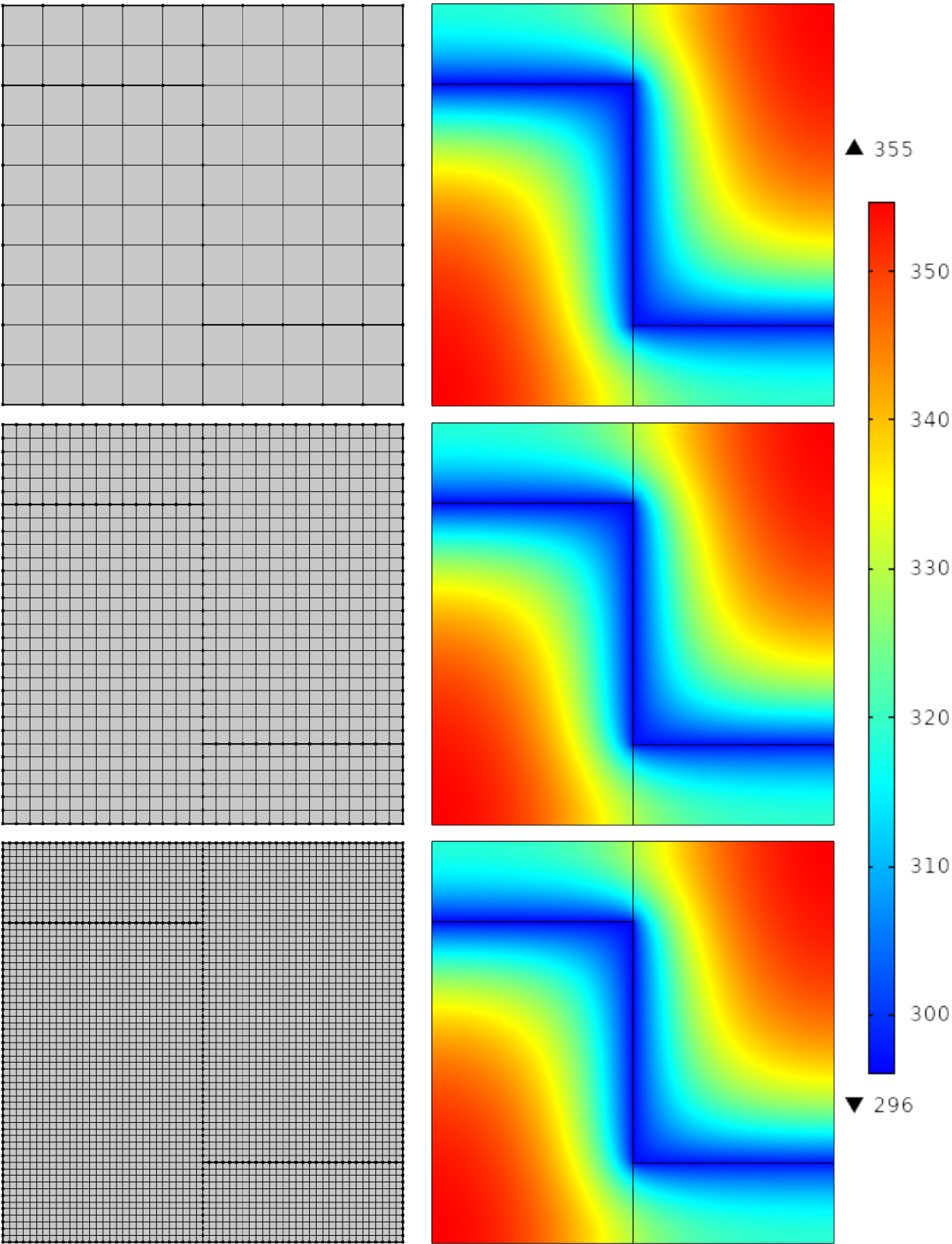}
    \caption{Stepped vasculature: $h$-refinement study. The figure shows a hierarchy of meshes on the left panel with corresponding temperature profiles on the right.}
    \label{Fig:App_Fig1_Stepped_vasculature}
\end{figure}

\begin{figure}
    \centering
    \includegraphics[scale=0.8]{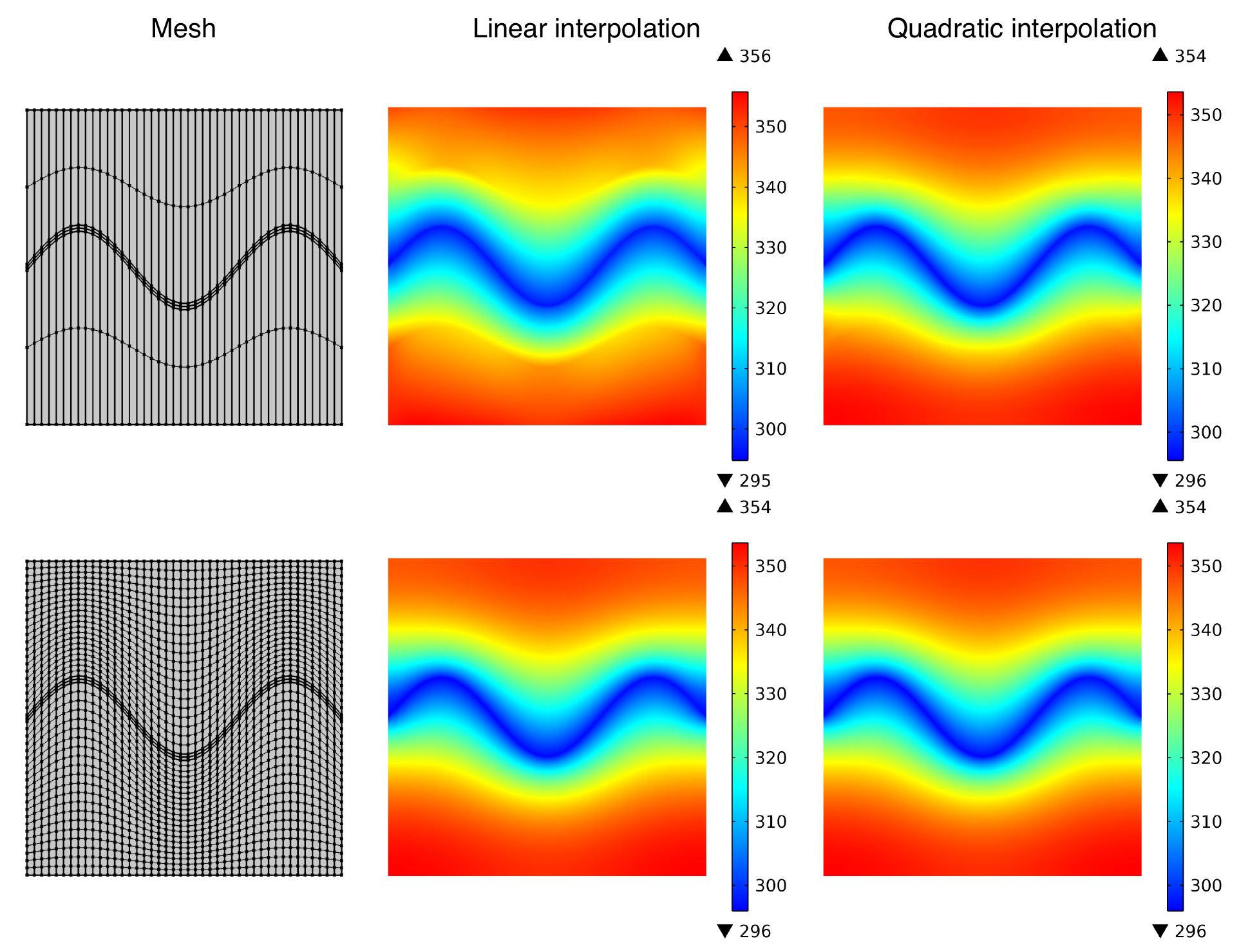}
    \caption{Sine-wave vasculature: $h$- and $p$-refinement study. The left panel shows two hierarchy meshes. The middle and right panels, respectively, show the corresponding temperature profiles under linear and quadratic interpolations.}
    \label{Fig:App_Fig2_Sine_wave}
\end{figure}

All the finite element simulations were conducted using the ``weak form" capability in COMSOL \citep{COMSOL}. The Galerkin finite element formulation was used; for details of this formulation, see \citep{nakshatrala2022qualitative}. To decide on the selection of elements and order of interpolation, we performed a numerical convergence study---varying the mesh size and the order of interpolation. The parameters were the same as in Table \ref{tab:params}. Figures \ref{Fig:App_Fig1_Stepped_vasculature} and \ref{Fig:App_Fig2_Sine_wave} show the convergence results. Quadrilateral elements with quadratic interpolation gave accurate results even on modestly refined meshes. Hence, we adopted the same element type and order of interpolation in generating all the finite element results presented in the main article.

\section*{DATA AVAILABILITY}
The data supporting this study's findings are available from the corresponding author upon reasonable request.

\section*{ACKNOWLEDGMENTS}
KBN acknowledges the support from the University of Houston through the High Priority Area Research Seed Grant. 
MKM acknowledges the support from the U.S. DOE Geothermal Technologies Office (Award \#: DE-EE-3.1.8.1). 
The authors' views and opinions expressed herein do not necessarily reflect those of the sponsors.

\bibliographystyle{unsrtnat}
\bibliography{Master_References}
\end{document}